\documentclass[a4paper,12pt]{article}
\usepackage[T1]{fontenc}
\usepackage[utf8]{inputenc}
\input{epsf.sty}
\usepackage{hyperref}
\usepackage[final]{epsfig}
\usepackage{color}
\usepackage{cite}
\usepackage{graphicx}
\usepackage{subcaption}
\usepackage{caption}
\usepackage{tikz}
\usepackage{makecell}
\usepackage{boldline}
\setcellgapes{3pt}
\usepackage{multirow}
\usepackage{mathtools}
\usepackage{adjustbox,lipsum}
\usetikzlibrary{decorations.pathreplacing}
\usepackage{amsmath, amsthm, amssymb, dsfont, mathrsfs}
\usepackage{comment}

\newtheorem {thm}{Theorem}[section]
\newtheorem {prop}[thm]{Proposition}
\newtheorem {lem}[thm]{Lemma}

\theoremstyle{definition}

\def\N{{\Bbb N}}

\def\Z{{\Bbb Z}}

\def\sm{\setminus}
\def\one{{\mathds 1}}

\newcommand{\C}{{\mathcal C}}

\newcommand{\e}{\mathrm e}
\renewcommand{\d}{\mathrm d}
\newcommand\restr[2]{{
  \left.\kern-\nulldelimiterspace 
  #1 
  \vphantom{\big|} 
  \right|_{#2} 
  }}

\def\0{{\bf 0}}

\def\dist{{\rm dist}}

\newcommand{\zd}{{{\mathbb Z}^d}}

\def\a{\alpha}

\def\phi{\varphi}

\def\g{\gamma}

\def\s{\sigma}

\def\o{\omega}

\def\D{\Delta}

\def\L{\Lambda}

\def\O{\Omega}

\def\T{\T}

\def\ot{\tilde{\omega}}

\def\ms{\text{ms}}

\def\is{T}

\def\tr{T}

\def\C{{\cal C}}

\def\PP{{\cal P}}

\def\FF{{\cal F}}

\def\EE{{\cal E}}
\def\QQ{{\cal Q}}
\def\WW{{\cal W}}

\begin{document}

\title{Gibbsianness of locally thinned random fields} 

\author{
Nils Engler\footnote{
Technische Universit\"at, Fakult\"at f\"ur Mathematik,  Stra\ss e des 17.~Juni 136, 10587 Berlin, Germany
\texttt{nils.engler@arcor.de}},
\, 
Benedikt Jahnel\footnote{
Weierstrass Institute for Applied Analysis and Stochastics, Mohrenstrasse 39
10117 Berlin, Germany
\texttt{Benedikt.Jahnel@wias-berlin.de}}
\, and
Christof K\" ulske\footnote{
Ruhr-Universit\"at Bochum, Fakult\"at f\"ur Mathematik, 44801 Bochum, Germany 
\texttt{Christof.Kuelske@ruhr-uni-bochum.de}}
}

\newcommand{\BJ}[1]{{\color{blue} #1}}
\newcommand{\NE}[1]{{\color{red} #1}}
\newcommand{\CK}[1]{{\color{green} #1}}

\maketitle

\begin{abstract}
We consider the locally thinned Bernoulli field on $\zd$, which is the lattice version of the Type-I Mat\'ern hardcore process in Euclidean space. It is given as the lattice field of occupation variables, obtained as image of an i.i.d.~Bernoulli lattice field with occupation probability $p$, under the map which removes all particles with neighbors, while keeping the isolated particles. 

We prove that the thinned measure has a Gibbsian representation and provide control on  its quasilocal dependence, both in the regime of small $p$, but also in the regime of large $p$, where the thinning transformation changes the Bernoulli measure drastically. Our methods rely on Dobrushin uniqueness criteria, disagreement percolation arguments~\cite{BeMa94}, and cluster expansions.   
\end{abstract}

\bigskip


\section{Introduction}\label{Sec_Intro}
Thinning transformations play a major role in the stochastic geometry of systems of 
point particles, see~\cite{rolski1991stochastic,moller2010thinning,isham1980dependent,last1993dependent,ball2005poisson,blaszczyszyn2019determinantal,bremaud1979optimal,jahnel2020probabilistic,baccelli2009stochastic}.
In that context a classical example is given when a point cloud is drawn according to a Poisson point process with homogeneous intensity in Euclidean space, from which 
afterwards all points are removed which have a neighbor at a distance less or equal than $1$, see~\cite{moller2010perfect,huber2009likelihood,matern1960spatial,matern2013spatial,penrose2001random,stoyan1988stochastic,andersen2016matern,baccelli2012extremal}.  

In this paper we consider a {\em discrete} version of such a transformation $T$ of {\em removal of non-isolates}, starting with the occupied sites drawn according to the i.i.d.~Bernoulli field on the integer lattice. The Bernoulli lattice field in itself and fine properties of the percolation transition driven by the occupation density $p$ is the object of a large literature, and ongoing research~\cite{MR1707339,MR2283880,MR3161674}. 

The Bernoulli lattice field also serves as a building block for more complex dependent processes of statistical physics and probability, which are derived from it. Let us mention bootstrap percolation (in which sites are added according to local growth  rules~\cite{adler1991bootstrap,MR3851832}), random walks on percolation clusters~\cite{MR3156983}, and diluted spin systems~\cite{MR1766342}.  
The latter two types of systems  are main subjects in the broader realm of disordered systems, see also~\cite{MR2252929,panchenko2013sherrington}. 

Specifically the {\em GriSing model} (where an Ising model is put on the occupied sites chosen as a realization of the Bernoulli lattice field) provides a somewhat surprising warning example of the appearance of a non-quasilocal measure, as the authors of~\cite{MR1766342} showed. This lack of quasilocality, {\em which is also termed non-Gibbsian behavior}, means that the system has finite-volume conditional probabilities with non-decaying dependence on variations of the boundary condition arbitrarily far away. For precise definitions see Section~\ref{Sec_Setting}. 
The non-quasilocality in the GriSing measure was shown to appear even in the regime of {\em subcritical} $p$, due to a mechanism related to Griffiths singularities~\cite{griffiths1969nonanalytic}
caused by arbitrarily large occupied clusters which appear at positive density. This shows that Gibbs properties and quasilocal dependence are subtle and may fail even in the absence of percolation.  

It was also discovered later that the {\em GriSing} measure is just example of the more general class of measures which may become non-quasilocal, namely the {\em joint measures of disordered spins systems} on the product space of disorder and spin variables~\cite{MR169647,kuhn1994critical, van2000comment}. Such systems may even possess full measure sets of discontinuity points for they specifications (that is their finite-volume measures in dependence on boundary conditions), 
which is a very strong form of singularity. This was shown in particular for the example of the joint measures of the random-field Ising model in the phase transition region on the lattice $\Z^3$ in~\cite{MR2060315}, building on~\cite{MR943702}.

For more studies of Gibbsian properties of transformed measures in probability and statistical physics under deterministic projection maps see~\cite{MR1012855,MR2097868,MR3858821,EnFeSo93,Fuz}. For related studies of Gibbs-non Gibbs transitions which are caused by stochastic dynamics, see~\cite{MR2185338,MR1889994,MR2643565,MR2976563,jahnel2017widom,kulske2007spin,MR4293773}. 

Let us come back to our Bernoulli lattice thinning process which we consider in the present paper. While the application of the thinning map $T$, as it {\em projects to isolates}, does not change the Bernoulli measure very much for small $p$, it changes the measure drastically for large $p$.  So one might conjecture that in particular the latter region is causing problems for a quasilocal Gibbsian description.  

As main result of the analysis of our paper we are however able to show that this is not the case, and in both regimes we have the regularity results of Theorem~\ref{theorem:main1} below, but for different reasons and with different proofs. Our proofs proceed via showing absence of phase transitions and regularity of the relevant internal systems (also known as first-layer measures under constraint), see the definition~\eqref{eq_2nd layer1}. For this we employ the suitably adopted mathematical-physics methods to prove uniqueness of infinite-volume measures which are {\em Dobrushin uniqueness criteria~\cite{Ge11}, disagreement percolation arguments~\cite{BeMa94}, and cluster expansions~\cite{friedli_velenik_2017}}. It turns out that there are some obstacles we need to overcome on the way to make this work, e.g., will it be necessary to go from a single-site description to a domino representation of the conditional first-layer measures, see Section~\ref{Sec_Proofs_Low}. We also provide a quantitative analysis 
and comparison of their effectiveness, in terms of numerical values for the regimes they can treat, 
see Section~\ref{Sec_Dis}.

Finally, let us note that our thinning map $T$, which is the {\em projection to isolates}, is accompanied by a natural companion map, namely the {\em projection to non-isolates} $T^*$, see Section~\ref{Sec_Setting}. Observe that the joint information of the images of both maps  provides a natural decomposition of the underlying i.i.d.~Bernoulli field. Since the latter Bernoulli field is trivially Gibbs as it even has no interaction, this suggests as a  first naive conjecture, that quasilocality of the projection map to {\em isolates} implies quasilocality of the projection map to {\em non-isolates}, too. We warn the reader that such a conclusion would be far too naive. On the contrary, our investigations in~\cite{jahnel2021gibbsianness} show that, indeed Gibbsianness of the projection to non-isolates fails for sufficiently large $p$. We highlight our findings in the following Table~\ref{table1}.

\makegapedcells
\begin{table}[!htpb]
  \centering
  \caption{Bernoulli $p$-projections: decomposition into isolates and non-isolates}
\label{table1}
\scalebox{1}{
  \begin{tabular}{cV{3}c|cV{3}c|c}
 & first-layer &  & Gibbs property of & \\
image measure & constraint model & range of $p$ & image measure& Reference\\
\hlineB{5}
$T\mu_p$&non-isolation &small& Gibbsianness&Thm.~\ref{theorem:main1} \\ \cline{3-5}
supported on& model on &  large & Gibbsianness&Thm.~\ref{theorem:main1} \\ \cline{3-5}
isolated sites& unfixed region&mid& Gibbs?& Sec.~\ref{Sec_Sim}\\ 
\hlineB{3}
$T^*\mu_p$&isolation &small& Gibbsianness&\cite[Thm.~2.2]{jahnel2021gibbsianness} \\ \cline{3-5}
supported on& model on &  large & non-Gibbsianness&\cite[Thm.~2.1]{jahnel2021gibbsianness} \\ \cline{3-5}
non-isolated sites& unfixed region&mid& sharp transition?& \\ 
\end{tabular}
}
\end{table}

The present paper is organized as follows. We present the setting and main results in Section~\ref{Sec_Setting}. In Section~\ref{Sec_Strategy} we present the strategy of the proofs. In Section~\ref{Sec_Supp} we elaborate on alternative strategies for parts of the proofs and evaluate their potential benefits for certain bounds in the parameter space. Here we also include a discussion on the intermediate regime that is not covered by our main results. Finally, in Section~\ref{Sec_Proofs} we present the proofs.

\section{Setting and main results}\label{Sec_Setting}
We consider the configuration space $\O = \{0,1\}^{\zd}$ equipped with the product topology and the associated Borel sigma-algebra $\FF$ for $d\ge 1$.  By $\mu_p \colon \FF \rightarrow \left[0,1\right]$ we denote the Bernoulli i.i.d.~product probability measure with density parameter $p \in \left[0, 1\right]$, i.e., $\mu_p(\o_i = 1) = p = 1 - \mu_p(\o_i = 0)$
independently for each $i \in \zd$. The event that $\o_i = 1$ is called {\em occupation} at $i \in \zd$, the complementary event is called {\em vacancy} at $i\in\Z^d$. We define the {\em isolation} event at site $i \in \zd$ by 
$$I_i := \{\o \in \O \colon \o_i = 1 \text{ and } \o_j = 0 \text{ for all } j \sim i\},$$ 
where $\sim$ denotes the usual neighborhood relation on $\Z^d$.

We further consider the associated deterministic thinning transformation $\is \colon \O \rightarrow \O$ given by
$$
(\is(\o))_i := \o'_i := \one\{\o\in I_i\},\qquad i \in \zd.
$$
In words, the transformation $\is$ removes all particles from the lattice, which have at least one neighboring particle. Note that $T$ is also a projection map since $T=T\circ T$. We further note that the complementary thinning $T^*(\o):=(\one\{\o\not\in I_i\})_{i\in \Z^d}$, which is also a projection, is considered in~\cite{jahnel2021gibbsianness}. Now, any $\o\in \O$ can be uniquely reconstructed from its joint images under the two maps as $\big((\one\{\o\in I_i\})_{i\in \Z^d},(\one\{\o\not\in I_i\})_{i\in \Z^d}\big)$. Next, let $\O' := \is(\O) \subset \O$ denote the space of particle configurations that obey the isolation constraint, and denote the image measure of $\mu_p$ under $T$ by 
$$
\mu'_p :=T\mu_p= \mu_p \circ \is^{-1}.
$$
Note that the mapping $\is$ defines a deterministic renormalization transformation in the sense of~\cite{EnFeSo93}, since it is local and maps translation-invariant measures onto translation-invariant measures.

\medskip
In this manuscript, we give answers to the question if the measure $\mu'_p$ is a Gibbs measure in the sense of existence of a quasilocal specification $\g'$ for $\mu'_p$. Recall that a {\em specification} $\g=(\g_\L)_{\L\Subset\Z^d}$ is a {\em consistent} and {\em proper} family of probability kernels, i.e., for all $\L\subset\D\Subset\Z^d$, $\o_\L\in \O_\L:=\{0,1\}^\L$ and $\hat\o\in \O$, we have that $\int_{\O}\g_\D(\d \tilde\o|\hat\o)\g_\L(\o_\L|\tilde\o)=\g_\D(\o_\L|\hat\o)$, and for all $\o_{\L^{\rm c}}\in \O_{\L^{\rm c}}$ we have $\g_\L(\o_{\L^{\rm c}}|\hat\o)=\one\{\o_{\L^{\rm c}}=\hat\o_{\L^{\rm c}}\}$ where $\hat\o_{\L^{\rm c}}$ denotes the restriction of $\hat\o$ to the volume $\L^{\rm c}$. A specification is called {\em quasilocal}, if for all volumes $\L\Subset\Z^d$ and local configurations $\hat\o_\L \in\O_\L$, the mapping $\o\mapsto\g_\L(\hat\o_\L|\o)$ is continuous with respect to the product topology on $\O$. We say that $\g$ is a specification for some random field $\mu$ on $\O$, if $\mu$ satisfies the DLR equations, i.e., for all $\L\Subset\Z^d$ and $\o_\L\in \O_\L$, we have that $\int_{\O}\mu(\d \tilde\o)\g_\L(\o_\L|\tilde\o)=\mu(\o_\L)$.
Here is our main result.
\begin{thm}[Gibbsianness for small and large $p$]\label{theorem:main1}
For all $d\ge 1$, there exist $0<p_1<p_2<1$ such that $\mu'_p$ is a Gibbs measure for $p\in [0,1]\setminus [p_1,p_2]$.
\end{thm}
The proofs for the different regimes require very different methods. We treat the small-$p$ case via cluster expansion and do not aim for explicit bounds on $p_1$. The cluster-expansion ansatz would also work for the large-$p$ case, however, this case, after some reformulations, can be treated via the less technical Dobrushin-uniqueness criterion. Using this, in particular, we can provide the following explicit lower bound on $p_2$.
\begin{prop}\label{Prop_DoBound}
Theorem~\ref{theorem:main1} holds for $p_2\le p^{\rm d}_{\rm c}(d)$, where
\begin{align*}
p^{\rm d}_{\rm c}(d)=\sup\{p\in(0,1)\colon 2(d-1)(d-2)(1-p^2)&+4(d-1)p(1-p)\\
&+2\frac{1-p}{1 - p(1-p)} + 6(d-1)(1-p)<1\}.
\end{align*}
\end{prop}
In the following section we give an overview of the strategies for the proofs. Note that, before we present the proofs in Section~\ref{Sec_Proofs}, we present some further results on the intermediate regime for $p$ and the optimality of the bound $p^{\rm d}_{\rm c}(d)$ in Section~\ref{Sec_Supp}.

\section{Strategy of proof}\label{Sec_Strategy}
The proofs depend on a two-layer approach. The {\em second-layer model} is given by $\mu_p'$, the thinned Bernoulli model with the hardcore constraint banning non-isolated sites as described above. Note that under the transformation $\is$, an occupied site in the thinned model determines its own value (occupied) and the values of all neighboring sites (unoccupied) of possible preimage configurations on the Bernoulli i.i.d.~field. Meanwhile, an unoccupied site after the thinning grants freedom in the choice of preimages in its neighborhood, as long as all occupied sites have at least one occupied neighbor. Given a thinned configuration, this observation allows to examine the i.i.d.~field only on the unfixed part of the lattice, where it is equipped with a hardcore non-isolation constraint. This is what we denote the {\em first-layer constraint model}. The main theorem, Theorem~\ref{theorem:main1}, is then proved in two steps. First, we construct regular versions of the conditional probabilities which are well-defined due to Gibbs-uniqueness in the first-layer constraint model. Here, the uniqueness can be guaranteed using Dobrushin-uniqueness bounds, disagreement-percolation thresholds and cluster-expansion techniques. Using this, it remains a small step to prove quasilocality of the constructed specification. 

\subsection{Transformations into first-layer constraint models}\label{Sec_1stlayer}
For any set of sites $\L\subset\Z^d$ we denote by $\L^c:=\Z^d\setminus \L$ its {\em complement} and by $\partial_-\L:=\{x\in \L\colon \text{there exists }y\in \L^c\text{ with }y\sim x\}$ its {\em inner boundary}. Moreover, we denote by $\L^o:=\L\setminus\partial_-\L$ the {\em interior} and by $\bar\L:=((\L^c)^o)^c$ the {\em extension} of $\L$. Finally, $\partial_+\L:=\bar\L\setminus \L$ denotes the {\em outer boundary} and $\partial\L:=\partial_-\L\cup \partial_+\L$ denotes the {\em thick boundary} of $\L$.

For any finite volume $\L\Subset\Z^d$, we wish to construct a candidate for the regular conditional probability $\g'_\L(\o'_{\L}| \o'_{\L^c})$ as a pointwise limit as $\D \uparrow \Z^d$ of the conditional probabilities of the transformed local configuration $\o'_\L$ given the event of some transformed annulus $\o'_{\D \sm \L}$ and some legitimate exterior non-transformed configuration $\o_{\D^c}$ beyond the annulus. Legitimacy here means that $\o'_{\D}\o_{\D^c} \in \tr^{-1}(\o')$. Here and in the sequel we will often make identifications of the form $\o=\o_\L\o_{\L^{\rm c}}$. Then, we define for $\L\subset\D$ and such boundary configurations, 
\begin{equation}\label{Eq0}
\begin{split}
\g'_{\o, \L}(\o'_\L|\o'_{\D\setminus \L})
:=&\, \frac{
\sum_{\o_\D}\mu_p(\o_\D)\one\{\tr_\D(\o_\D \o_{\D^c})=\o'_\D\}}{
\sum_{\o_{\D\setminus\L^o}}\mu_p(\o_{\D\setminus\L^o})\one\{\tr_{\D \setminus \L}(\o_{\D\setminus\L^o} \o_{\D^c})=\o'_{\D\setminus\L}\}
}\\
=&\, \frac{
\sum_{\o_{\D\setminus\L^o}}\mu_p(\o_{\D\setminus\L^o})\one\{\tr_{\D \setminus \L}(\o_{\D\setminus\L^o} \o_{\D^c})=\o'_{\D\setminus\L}\}
F[\o'_\L](\o_{\partial \L}) }{
\sum_{\o_{\D\setminus\L^o}}\mu_p(\o_{\D\setminus\L^o})\one\{\tr_{\D \setminus \L}(\o_{\D\setminus\L^o} \o_{\D^c})=\o'_{\D\setminus\L}\}
},
\end{split}
\end{equation}
where $\mu_p(\o_{\L})=\prod_{i\in\L}p^{\o_i}(1-p)^{1-\o_i}$ is the Bernoulli measure in the volume $\L$, we wrote $\tr_\L(\o)$ instead of $(\tr(\o))_\L$, and 
$$
F[\o'_\L](\o_{\partial \L}):= \sum_{\o_{\L^o}}\mu_p(\o_{\L^o})\one\{\tr_{\L}(\o_{\L^o}\o_{\partial \L} )=\o'_{ \L} \}
$$
is a local function. Let us recall the following general result about the specification property, whose proof is based on martingale-convergence arguments.
\begin{lem}[{\cite[Lemma 3.3]{jahnel2021gibbsianness}}\, ]\label{lem_specification}
Assume that, given $\L \Subset \Z^d$ and $\o' \in \O'$, we have that the limit $\g'(\o'_\L | \o'_{\L^c}) := \lim_{\D \uparrow \zd} \g'_{\o, \D}(\o'_\L|\o'_{\D\sm \L})$ exists and is independent of $\o \in T^{-1}(\o')$. Then, $\g'$ is a specification for $\mu'_p$.
\end{lem}
Hence, we need to guarantee the existence of a limiting object $\g'_\L(\o'_{\L}| \o'_{\L^c})$ that is independent of the external boundary condition $\o$. For this, our strategy is to invoke Gibbs-uniqueness criteria. In order to do this, first note that we can uniquely identify $\o'$ with the subset of its occupied sites in $\Z^d$. With a slide abuse of notation, we can then see that the extension $\bar\o':=\overline{\o'}$ of $\o'$ is a {\em fixed region} for the first-layer constraint model as defined below, in the sense that, under the transformation, there is no choice for the Bernoulli field in how to realize $\o'$. 

In view of this, we consider a general choice $S\subset\Z^d$ for the {\em unfixed region} (the complement of $\bar\o'$) and introduce the following specification associated to the first-layer constraint model on $\{0,1\}^S$,
\begin{equation}\label{eq_2nd layer1}
\g^S_{\D}(\o_\D|\o_{\D^c})
:=\frac{
\mu_p(\o_{\D\cap S})\one\{\o_{\D\cap S}\o_{\D^c\cap S}\text{ is $T$-feasible on } \D\cap S\}}{
\sum_{\ot_{\D\cap S}}\mu_p(\ot_{\D\cap S})\one\{\ot_{\D\cap S}\o_{\D^c\cap S}\text{ is $T$-feasible on } \D\cap S \}},\qquad\D\subset\Z^d.
\end{equation}
Here, a configuration $\o\in \O$ is called {\em $\is$-feasible} on a set $\D\cap S$ if all occupied sites of $\o$ in $\D\cap S$ have at least one occupied neighbor, which may lie in $\bar\D\cap S$. In particular, with this definition, 
\begin{equation*}
\g'_{\o, \L}(\o'_\L|\o'_{\D\setminus \L})=\g^{(\bar\o')^{\rm c}}_\D(F[\o'_\L]|\o_{\D^c}),\qquad \D\Subset\Z^d,
\end{equation*}
where we used that in the fixed area we see cancellations. 
Then, we have the following propositions that we prove in Section~\ref{Sec_Proofs}. 
\begin{prop}[Low-density Gibbsianness]\label{proposition_low}
There exist $0<p_1<1$ such that for all $0\leq p\leq p_1$ and $\o'\in \Omega'$ the limit $\lim_{\D\uparrow\Z^d}\g^{(\bar\o')^c}_\D(F[\o'_\L]|\o_{\D^c})=:\g'(\o'_\L|\o'_{\L^c})$ exists independently of $\o\in T^{-1}(\o')$. Moreover, $\g'$ is a quasilocal specification for $\mu'_p$. 
\end{prop}
The proof of Proposition~\ref{proposition_low} is based on cluster-expansion techniques as the specification kernel of the first-layer constraint model fails to satisfy Dobrushin's condition of weak dependence, due to the non-isolation constraint, and will be presented in Section~\ref{Sec_Proofs_Low}.
\begin{prop}[High-density Gibbsianness]\label{proposition_high}
There exist $0<p_2<1$ such that for all $p_2\leq p\leq 1$ and $\o'\in \Omega'$ the limit $\lim_{\D\uparrow\Z^d}\g^{(\bar\o')^c}_\D(F[\o'_\L]|\o_{\D^c})=:\g'(\o'_\L|\o'_{\L^c})$ exists independently of $\o\in T^{-1}(\o')$. Moreover, $\g'$ is a quasilocal specification for $\mu'_p$. 
\end{prop}
The proof of Proposition~\ref{proposition_high} is based on Dobrushin-uniqueness techniques and will be presented in Section~\ref{Sec_Proofs_High}. 
Before we exhibit the proofs, in the following Section~\ref{Sec_Supp}, we present some supplementary results.

\section{Alternative bounds and intermediate regimes}\label{Sec_Supp}
In this section we present further results on the bounds for $p_2$ as well as on the behavior of the system for intermediate values of $p$. 

\subsection{Disagreement-percolation bounds}\label{Sec_Dis}
\begin{figure}[t]
\centering
\includegraphics[width=1\textwidth]{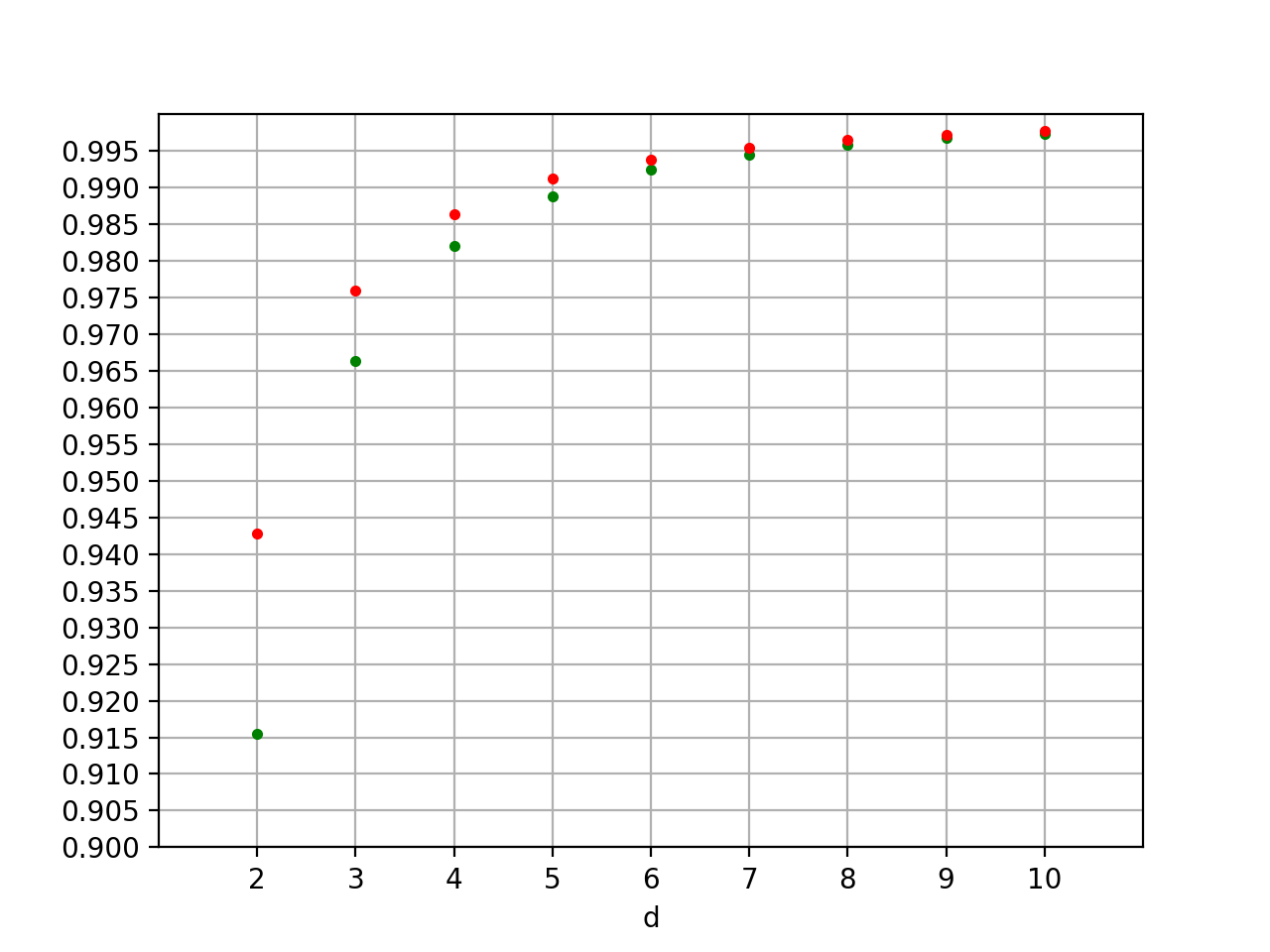}\hfill
\caption{Comparison of $d\mapsto p^{\rm d}_{\rm c}(d)$ (green) and $d\mapsto p^{\rm p}_{\rm c}(d)$ (red).}
\label{fig:figure_dob_thresh}
\end{figure}

The lower bound $p^{\rm d}_{\rm c}(d)$ of Proposition~\ref{Prop_DoBound} for the high-density Gibbsian regime, is a consequence of the Dobrushin-uniquness criterion for the first-layer constraint model~\eqref{eq_2nd layer1}. It guarantees unique existence of the infinite-volume Gibbs measure for~\eqref{eq_2nd layer1}, uniformly over the unfixed area $S$. However, there are alternative approaches in order to establish the unique existence of this infinite-volume first-layer constraint model, e.g., disagreement-percolation criteria. Let us next present a corresponding bound and discuss the relation to the Dobrushin-uniqueness bound. 
\begin{prop}\label{Prop_DisBound}
Let $d\geq 2$, then, for $p>p^{\rm p}_{\rm c}(d)$, with
$$
p^{\rm p}_{\rm c}(d)=\sqrt{\frac{2d^2 + 2d - 4 }{2d^2 + 2d - 3}},
$$
the first-layer constraint model $\g^S$, as defined in~\eqref{eq_2nd layer1}, admits a unique infinite-volume Gibbs measure, for all unfixed areas $S$. 
\end{prop}
Let us note that $p^{\rm p}_{\rm c}(2)=\sqrt{8}/3 \approx 0.9428$ and $p^{\rm p}_{\rm c}(3)=\sqrt{20}/\sqrt{21} \approx 0.9759$. On the other hand, $p^{\rm d}_{\rm c}(2)\approx 0.9155$ and $p^{\rm d}_{\rm c}(3)\approx 0.9663$ and this trend, that the Dobrushin criterion provides better bounds with a decreasing difference, as the dimension grows, can also be observed by further simulations, see Figure~\ref{fig:figure_dob_thresh}.
However, let us note that the bound used to derive $p^{\rm p}_{\rm c}$, is also certainly not optimal since it is based on a general criterion~\eqref{eq_percolation_bound} for percolation via maximal graph degrees, see Section~\ref{Sec_Prop_DisBound} for details. Indeed, incorporating Monte--Carlo simulation results from~\cite[Table 1]{malarz2005square} for critical values for percolation in the interaction graph of the dimer representation of our model, lead for example to the smaller value $p^{\rm p}_{\rm c}(2)=\sqrt{8}/3 \approx 0.8438$. We present the proof of Proposition~\ref{Prop_DisBound} in Section~\ref{Sec_Prop_DisBound}.

\subsection{Some computations for the intermediate regime}\label{Sec_Sim}
An open question is whether Gibbsianness holds in the intermediate regime, that is, e.g., for $p \approx 1/2$. The standard way of proving non-Gibbsianness is to determine points of essential discontinuity, i.e., certain configurations $\o'$, such that alterations of spins arbitrarily far away from $\L$ change the value of $\g'_\L(\o'_\L|\o'_{\L^c})$ by an amount greater than some fixed $\epsilon > 0$. We are hence looking for thinned configurations, such that the Bernoulli i.i.d.~measure with non-isolation constraint on the unfixed part $(\bar\o')^{\rm c}$ is most likely to exhibit a phase transition. One special candidate is given by the checkerboard (or alternating) configuration $\o'_{\rm{alt}}$, where each site $(\o'_{\rm{alt}})_i$ is occupied if and only if $\sum_{k=1}^d i_k$ is odd. Note that we have $\o'_{\rm{alt}} \in \O'$. However, since $(\bar\o'_{\rm{alt}})^{\rm c} = \emptyset$, there cannot be transport of information through the annulus $\D \sm \L$ in this case since there are no internal spins allowing for a phase transition to occur. We note that such transport of information gets more likely, the larger the unfixed area becomes. Hence, it seems reasonable to study the completely unoccupied configuration $\o'_{\text{zero}} \in \O'$ with $(\bar\o'_{\text{zero}})^{\rm c}= \zd$ as a potential point of essential discontinuity.

\begin{figure}[!htpb]
\centering
\includegraphics[width=.333\textwidth]{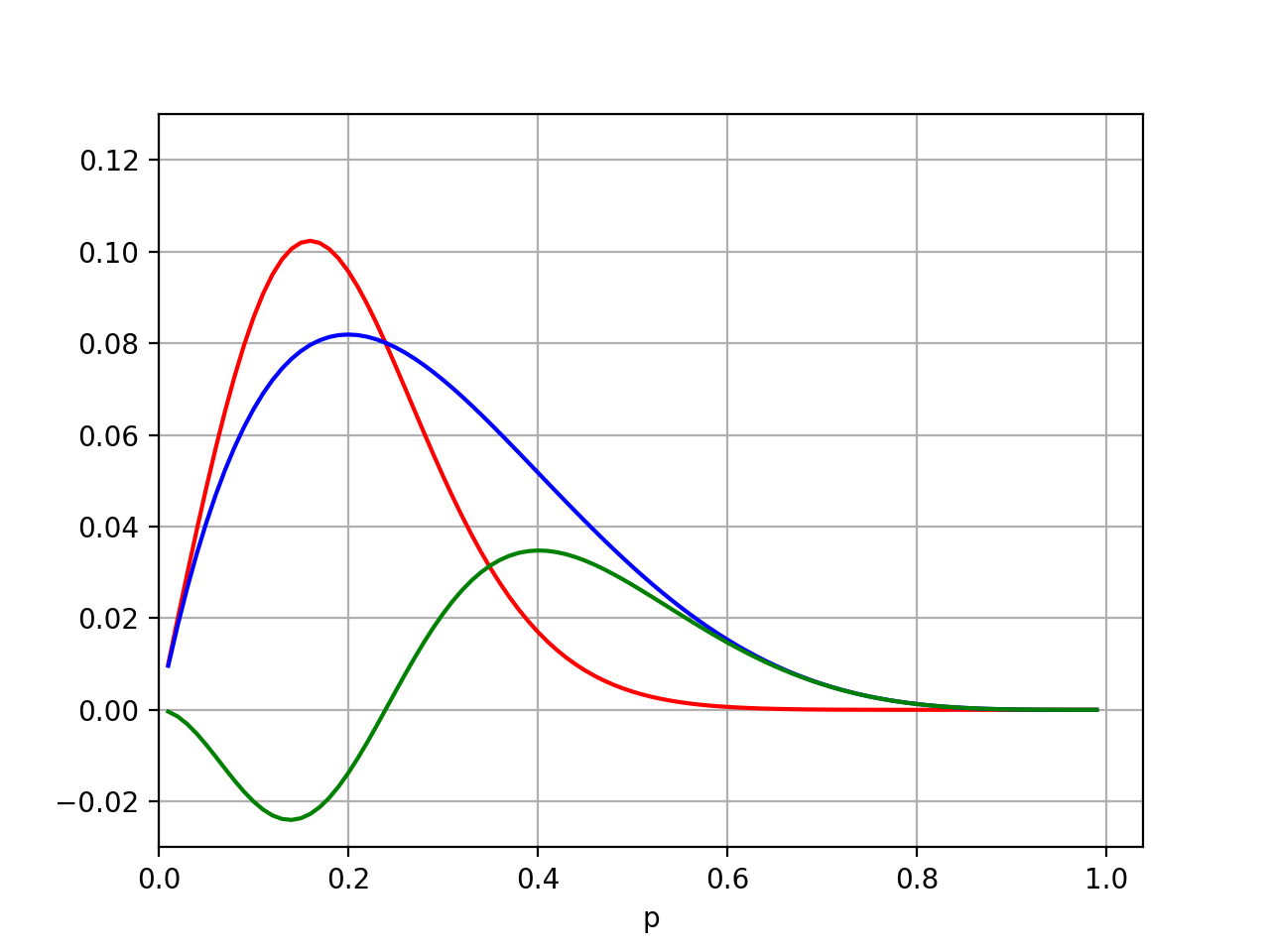}\hfill
\includegraphics[width=.333\textwidth]{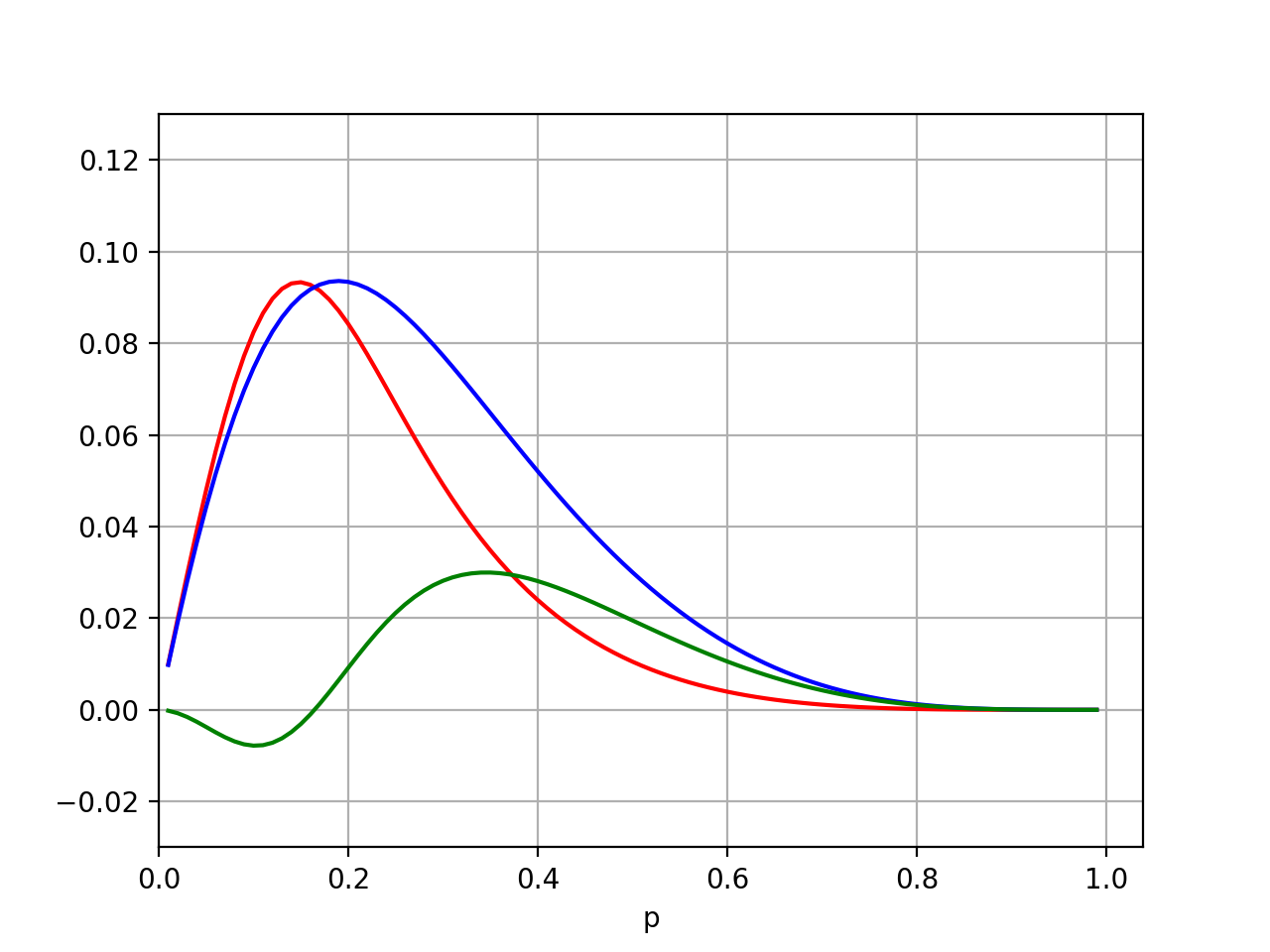}\hfill
\includegraphics[width=.333\textwidth]{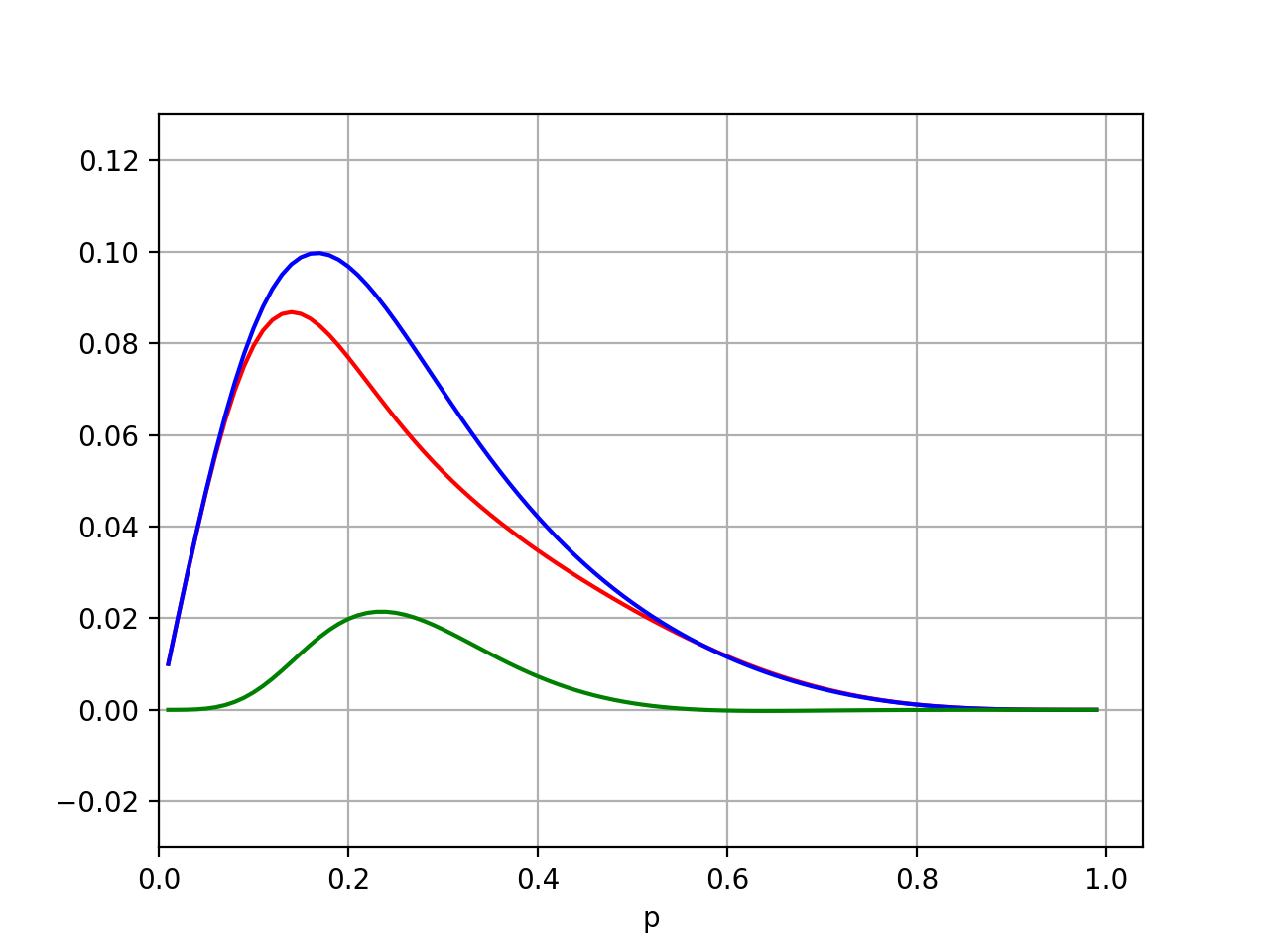}
\caption{Computation of the conditional probability in (\ref{Eq0}) for $\protect d=2$ with $\protect \L = \{0\}$, $\protect \D = B_3, B_4, B_5$ (cubes around the origin with sidelength $i=3,4,5$), occupied origin, and unoccupied surrounding on $\D \setminus \L$. The red lines correspond to the unoccupied boundary condition $\o$ on $\D^c$, the blue lines to the occupied boundary, while the green lines represent the difference of the two. It can be observed that the difference of the conditional probabilities of different boundary conditions appears to decrease uniformly in $\protect p$ with increasing $\protect \D$.}
\label{fig:figure3}
\end{figure}

For this purpose, we have written a script in order to compute the exact values of the conditional probabilities in (\ref{Eq0}) in dimension $d=2$ for $\L = \{0\}$ (origin), $\D = B_k, k=3,4,5$ (cubes around the origin of side length $k$), a fully unoccupied annulus $\o'_{\D \sm \L} = (\o'_{\text{zero}})_{\D \sm \L}$ and either a fully occupied, or fully unoccupied boundary condition off the annulus. We present the resulting exact calculations in Figure \ref{fig:figure3} and note that they seem to suggest that, at least in two dimensions, there is no phase transition even in the intermediate density regime. The code for the computations can be found at \cite{boxes}.
Let us finally note that, in order to prove or disprove that the Gibbs property persists for all $p$, different methods have to be developed and more research is necessary. 

\section{Proofs}\label{Sec_Proofs}
We will often suppress the dependence on $p$ in the remainder of the paper, whenever there is no risk for ambiguity. 

\subsection{Proof of Proposition~\ref{proposition_low}}
\label{Sec_Proofs_Low}
Let $S\subset\Z^d$ be an unfixed area that will be fixed for the most part of the section. We consider subsets $\L\subset\D\Subset S$. In order to ease notation in the remainder of this section, any operation such as $\L^{\rm c}$, or $\partial_-\L$ should be understood with respect to $S$. For example $\L^{\rm c}=\{x\in S\colon x\notin \L\}$ or $\partial_-\L=\{x\in \L\colon \text{there exists }y\in \L^c\text{ with }y\sim x\}$. 
\subsubsection{Cluster expansion}
The proof proceeds via cluster expansion on the annulus $\D\setminus \bar\L$. 
Assuming $\D$ to be sufficiently large, we can split the outer boundary of $\D\setminus \bar\L$ into an inner and an outer part, i.e., for any $\o\in \O$, $\o_{\partial_+(\D\setminus\bar\L)}=\o_{\partial_+\L}\o_{\partial_+\D}$. 
We like to stress that boundary sites are never in $\Z^d\setminus S$. The idea is to derive an expansion for the {\em partition function}
\begin{equation}\label{eq_partition_cluster_proof}
Z_{\D\setminus \bar\L}(\o_{\partial\L}\o_{\partial_+\D}) := \sum_{\o_{\D\setminus \bar\L}}\mu(\o_{\D\setminus \bar\L})\one\{\o_{\partial\L}\o_{\D\setminus \bar\L}\o_{\partial_+\D} \text{ is $T$-feasible on } \D\setminus\L\},
\end{equation}
where we want to highlight the fact that the inner and the outer part of the boundary are treated differently. The reason for this is that the quantities in (\ref{Eq0}) and (\ref{eq_2nd layer1}), whose limit we would like to investigate, require feasibility only on $\D$. By taking advantage of cancellations, we then show that for any two $\o_{\partial\L}, \tilde \o_{\partial\L}$
\begin{equation}\label{Goal_1}
\begin{split}
\lim_{\D \uparrow \Z^d}\frac{Z_{\D\setminus \bar\L}(\o_{\partial\L} \o_{\partial_+\D})}{Z_{\D\setminus \bar\L}(\tilde \o_{\partial\L} \o_{\partial_+\D})}
\end{split}
\end{equation}
exists and is independent of $\o_{\partial_+\D}$. 

\paragraph{Setting up the cluster expansion:}
To start, we define the set $J_i(\o_{\partial\L} \o_{\partial_+\D})$ of configurations in $\D\setminus \bar\L$, such that the site $i\in \D\setminus \bar\L$ has an isolated occupied neighbor in $\D\setminus\L$, i.e., 
$$
J_i(\o_{\partial\L} \o_{\partial_+\D}) := \{\o_{\D\setminus \bar\L} \colon \text{there exists } j \sim i, j \in \D\setminus\L, \text{ such that } \o_j = 1\text{ and } \o_k = 0 \text{ for all } k \sim j \}.
$$
In particular, $\o_{\D\setminus \bar\L}\in J_i(\o_{\partial\L} \o_{\partial_+\D})$ implies $\o_i = 0$. The partition function (\ref{eq_partition_cluster_proof}) now becomes
\[
\sum_{\o_{\D\setminus \bar\L}}\mu(\o_{\D\setminus \bar\L})\prod_{i \in \D\setminus \bar\L} (1 - \one\{\o_{\D\setminus \bar\L}\in J_i(\o_{\partial\L} \o_{\partial_+\D})\}),
\]
where the product replaces the constraint that each site is not allowed to have an isolated neighbor. Then, we can rewrite
$$
\prod_{i\in \D\setminus \bar\L}(1 - \one\{\o_{\D\setminus \bar\L}\in J_i(\o_{\partial\L} \o_{\partial_+\D})\}) = \sum_{W \subset \D\setminus \bar\L}\prod_{i\in W}(- \one\{\o_{\D\setminus \bar\L}\in J_i(\o_{\partial\L} \o_{\partial_+\D})\}). 
$$
Let us define a notion of distance on $S$. We define ${\rm d}_S(i,j)$ to be the length of the shortest path in $S$, which starts in $i$ and ends in $j$. We then denote by $B^S_n(i)$ the associated ball of radius $n$ centered at $i\in S$. 
Using this, we can decompose each subset $W$ into its maximally connected components $W_1, \dots, W_n$ with respect to the graph on $\D\setminus \bar\L$ in which $i, j \in W$ are {\em connected} if and only if $i\in B_4^S(j)$. For this, we define the {\em dependence set} of $W_i$ to be $\overline W_i := \bigcup_{j \in W_i}B^S_2(j)$. We call $W_i, W_j$ {\em compatible} if and only if  $\overline W_i \cap \overline W_j = \emptyset$. The sets $W_1, \dots, W_n$ play the role of polymers with a hardcore interaction given by compatibility. 

\medskip
Due to the construction, the random variables $\prod_{i\in W_k}(-\one_{J_i(\o_{\partial\L} \o_{\partial_+\D})})$, $k = 1, \dots, n$ are independent with respect to $\mu$, and hence
\begin{equation*}\begin{split}
Z_{\D\setminus \bar\L}(\o_{\partial\L} \o_{\partial_+\D})& =
\sum_{\substack{W_1, \dots, W_n \\ \text{pw.~comp.}}}\prod_{k=1}^n \sum_{\o_{\D\setminus \bar\L}}\mu(\o_{\D\setminus \bar\L})\prod_{i\in W_k}(-\one\{\o_{\D\setminus \bar\L}\in J_i(\o_{\partial\L} \o_{\partial_+\D})\})\\
&=:\sum_{\substack{W_1, \dots, W_n \\ \text{pw.~comp.}}}\prod_{k=1}^n z^p_{W_k}(\o_{\partial\L} \o_{\partial_+\D}),
\end{split}
\end{equation*}
where the first sum runs over all possible families of pairwise-compatible subsets of $\D\setminus\bar\L$. 

\medskip
Next, we define the set of {\em polymers} by 
\begin{equation*}
\begin{split}
\Gamma_{\D\setminus \bar\L}(\o_{\partial\L} \o_{\partial_+\D}):= \{W_k& \subset \D\setminus\bar\L\colon \text{ there exists } W \subset   \D\setminus\bar\L  \text{ and } \o_{\D\setminus \bar\L} \text{ such that }W_k \text{ is a maximally}\\ 
& \text{connected component of } W\text{ and }  \prod_{i\in W}\one\{\o_{\D\setminus \bar\L}\in J_i(\o_{\partial\L} \o_{\partial_+\D})\}=1\},
\end{split}
\end{equation*} 
where the connectedness is in the above sense. 
We need to treat polymers whose dependence sets intersect $\partial\L$ separately and thus define 
$$\QQ_{\D\setminus \bar\L}(\o_{\partial\L} \o_{\partial_+\D}) := \{Q \in \Gamma_{\D\setminus \bar\L}(\o_{\partial\L} \o_{\partial_+\D}) \colon \overline Q \cap \partial\L \neq \emptyset\},$$
the set of such intersecting polymers. The reason for doing this is that these polymers are not exponentially suppressed in their full volume, which is why their corresponding cluster expansion does not necessarily converge, see the estimates around the display~\eqref{eq_Q_k} below. Now, in order to maintain compatibility of all polymers, we define for a given collection of pairwise-compatible polymers $Q_1,\dots, Q_n \in \QQ_{\D\setminus \bar\L}(\o_{\partial\L}\o_{\partial_+\D})$ the set 
$$\WW_{Q_1,\dots, Q_n, \D\setminus \bar\L}(\o_{\partial_+\D}) := \{W \in \WW_{\D\setminus \bar\L}(\o_{\partial_+\D}) \colon \overline W \cap \overline Q_i = \emptyset \text{ for all } i= 1,\dots, n\},$$ 
of polymers compatible with that collection, where $\WW_{\D\setminus \bar\L}(\o_{\partial_+\D}) := \Gamma_{\D\setminus \bar\L}(\o_{\partial\L} \o_{\partial_+\D}) \sm \QQ_{\D\setminus \bar\L}(\o_{\partial\L} \o_{\partial_+\D})$. 
Then, we can write
\begin{equation*}\begin{split}
Z_{\D\setminus \bar\L}(\o_{\partial\L} \o_{\partial_+\D}) 
=  &\sum_{\substack{Q_1, \dots, Q_n \in \QQ_{\D\setminus \bar\L}(\o_{\partial\L}\o_{\partial_+\D}) \\ \text{pw.~comp.}}}\prod_{k=1}^n z^p_{Q_k}(\o_{\partial\L} \o_{\partial_+\D})
\sum_{\substack{W_1, \dots, W_m \in \WW_{Q_1,\dots, Q_n,\D\setminus \bar\L}(\o_{\partial_+\D}) \\ \text{pw.~comp.}}}\prod_{j=1}^m z^p_{W_j}( \o_{\partial_+\D})
\end{split}
\end{equation*}
and derive a convergent cluster representation for the polymers in $\WW_{Q_1,\dots, Q_n,\D\setminus \bar\L}(\o_{\partial_+\D})$ for sufficiently small $p$.

\paragraph{Convergence of the cluster expansion:} 
Let us suppress the dependence on $Q_1,\dots, Q_n$ for notational convenience in this part.
In order to derive the cluster representation, let us define for any {\em cluster} $C =\{W_1,\dots, W_n\}^{\ms}$, i.e., a multiset of pairwise-compatible polymers in $\WW_{\D\setminus \bar\L}(\o_{\partial_+\D})$, and boundary condition $\o_{\partial_+\D}$, the {\em cluster potential}
\begin{equation}\label{Eq_clus2}
\Phi_{\D\setminus \bar\L, \o_{\partial_+\D}}^p(C):= \left(\prod_{W \in \WW_{\D\setminus \bar\L}(\o_{\partial_+\D})}\frac{1}{n_C(W)!}\right)
\left(\sum_{\substack{G \subset G_n \\ \text{connected}}}
\prod_{\{i, j\} \in G}\zeta(W_i,W_j) \right)
\left(\prod_{k = 1}^{n}z^p_{W_k}(\o_{\partial_+\D}) \right),
\end{equation}
where, $n_C \colon \Gamma \to \mathbb{N}_{0}$ denotes a map that assigns to each polymer the number of occurrences in the cluster $C$. The function $\zeta(W_i,W_j)$ equals $0$ if $W_i$ and $W_j$ are compatible and $-1$ otherwise, while the sum is over all connected subgraphs of the complete graph $G_n = (V_n, E_n)$ on $n$ vertices. Here, the notion subgraph refers to a graph $G = (V, \EE)$, for which $V = V_n$ and $\EE \subset E_n$. 

\begin{figure}[t]
\centering
\includegraphics[width=\textwidth]{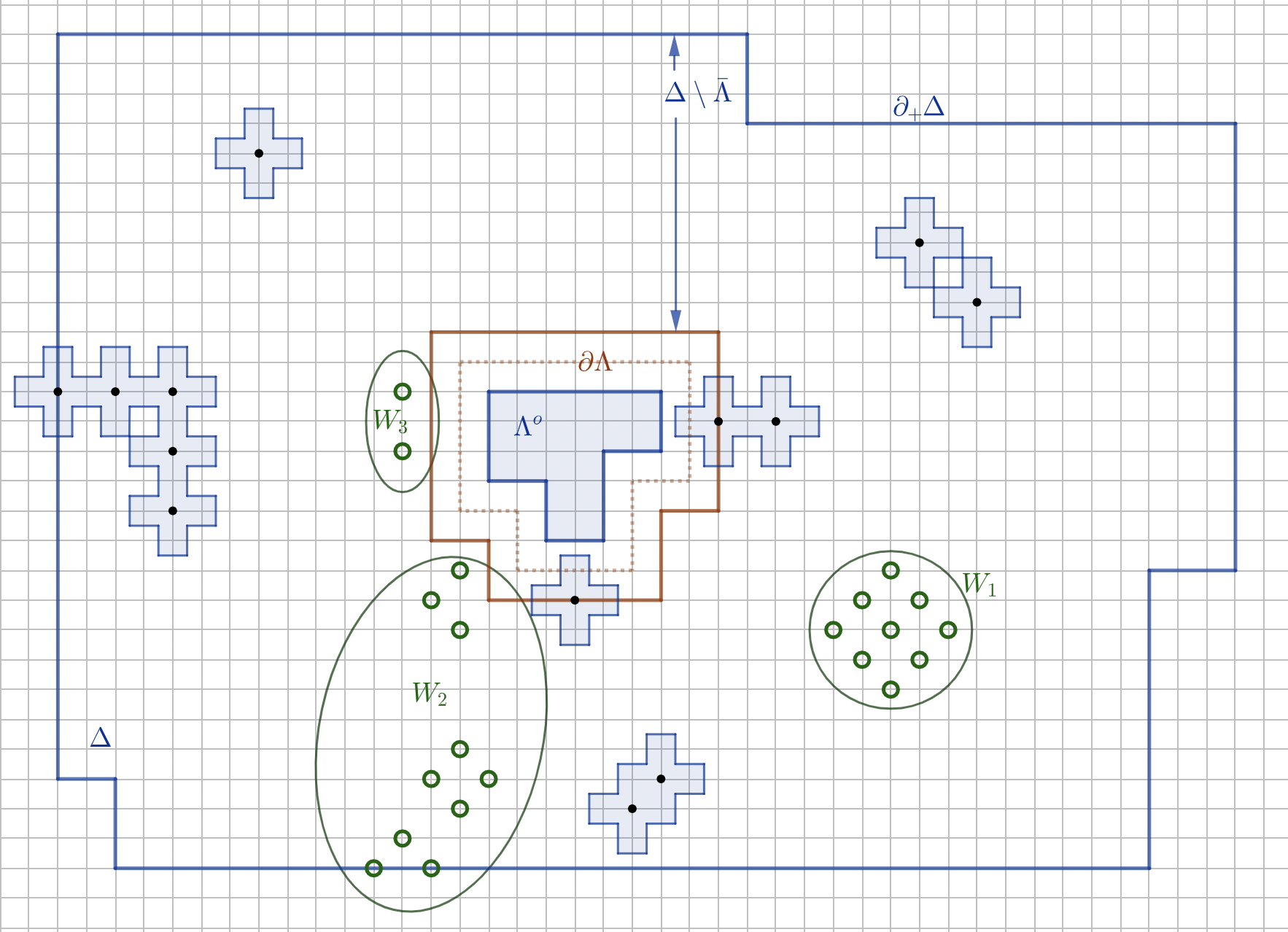}\hfill
\caption{Given thinned configurations $\protect \o'_{\L}, \o'_{\D \sm \L}$ and unthinned outer boundary configuration $\o_{\D^c}$, the fixed area $\bar\o'_{\D \setminus \L}$ is given by the crosses hatched in blue. There are three polymers (green). $\protect W_3$ at the top is due to isolations on the set $\partial_+\L$.}
\label{fig:figure1}
\end{figure}

We want to employ the criterion~\cite[Theorem 5.4.]{friedli_velenik_2017} in order to establish convergence of the cluster potentials. For this, the main ingredient is the following estimate for the polymer weights $W_i \in \WW_{\D\setminus \bar\L}(\o_{\partial_+\D})$,
\[
|z^p_{W_k}(\o_{\partial_+\D})| \leq p^{|L_{W_k}|} \leq p^{|W_k|/(2d)},
\]
where
 $$L_{W_k} := \{j \in \overline W_k\colon \text{ for all } i \sim j \text{ either } i \in W_k \text{ or } i \in (\D\setminus\bar\L)^{\rm c} \text{ with } \o_i = 0\}$$
denotes the set of sites completely surrounded by $W_i$ or by zeros on the boundary. As an example, consider the four sites enclosed by the polymer $W_1$ in Figure \ref{fig:figure1}.
We verify the condition of~\cite[Theorem 5.4.]{friedli_velenik_2017} for the volume function $a(W)= |W|$. Indeed, for any  polymer $W^* \in \WW_{\D\setminus \bar\L}(\o_{\partial_+\D})$, we find that
\begin{equation}\label{eq_s}
\begin{split}
\sum_{\substack{W \in \WW_{\D\setminus \bar\L}(\o_{\partial_+\D}): \\ \overline W \cap \overline W^* \neq \emptyset}}z^p_{W}(\o_{\partial_+\D})\e^{|W|}&\leq \sum_{\substack{W \in \WW_{\D\setminus \bar\L}(\o_{\partial_+\D}): \\ \overline W \cap \overline W^* \neq \emptyset}}p^{|W|/(2d)}\e^{|W|}\\
&\leq  |\overline W^*| \sum_{k \geq 1} \e^{k(1+\log(p)/2d)}  |\{W : 0 \in \overline W, |W| = k\}|\\
&\leq ((2d)^2+1)|W^*| \sum_{k \geq 1} \e^{k(1+\log(p)/2d))} (9^{d})^{(k+1){9^d}}\\
&= ((2d)^2+1)(9^d)^{9^d}|W^*| \sum_{k \geq 1} (\e^{1+\log(p)/2d} (9^d)^{9^d})^{k}
\leq  |W^*|
\end{split}
\end{equation}
for all $p \leq q_1$, with $0<q_1$ chosen sufficiently small. In the third inequality, we have bounded the number of polymers $W$ of size $k$ that contain $o$ in its dependence set by the number of connected graphs $H = (W \cup \{0\}, L)$ on $k+1$ nodes, where $\{i, j\} \in L$ if and only if $j\in B_4(i)$, where $B_n=B^{\Z^d}_n$ like above. We note that the possible number of incident edges for each node is bounded by $|B_4(o)| \leq 9^d$. Since each edge is shared by two nodes, we can bound the number of edges $|L|$ by $(k+1) 9^d/2$. Moreover, since for each connected graph and each starting node, there exists a path visiting each edge exactly twice, we can bound the number of graphs by $(9^d)^{2|L|}$.

\medskip
At this point we also see that there is monotonicity in the sense that when removing sites from the grid, i.e., if we consider a subset $\tilde S \subset S$, the number of polymers of a given size that are also incompatible with $W^*$ decreases, which is why convergence is ensured for $p\le q_1$ uniformly in $S$.

\medskip
Finally, an application of~\cite[Theorem 5.4.]{friedli_velenik_2017} ensures that the representation 
\begin{equation}\label{Eq_clus}\begin{split}
Z_{\D\setminus \bar\L}(\o_{\partial\L} \o_{\partial_+\D})=\sum_{\substack{Q_1,\dots, Q_n \in \QQ(\o_{\partial\L} \o_{\partial_+\D}) \\ \text{pw.~comp.}}}\prod_{k=1}^n z^p_{Q_k}(\o_{\partial\L} \o_{\partial_+\D})
\exp \Big( \sum_{\substack{C \in P(\WW_{Q_1,\dots, Q_n, \D\setminus \bar\L}(\o_{\partial_+\D}))}}\Phi_{\D\setminus \bar\L,\o_{\partial_+\D}}^p(C)\Big)
\end{split}
\end{equation}
is well defined, where $P(\WW)$ denotes the set of all clusters in $\WW$. 

\paragraph{Suppression of large clusters:}
We now work towards our goal~\eqref{Goal_1}, the convergence of the fraction of partition functions for different interior boundary conditions. For this purpose, we wish to bound the contribution of large clusters in the expansion \eqref{Eq_clus}. Let $\overline C=\bigcup_{i=1}^n W_i$ denote the {\em support} of the cluster $C =\{W_1,\dots, W_n\}^{\ms}$, fix $x\in \D\setminus \bar\L$ and write $\phi(C)$ for the first two combinatorial factors in (\ref{Eq_clus2}). We then get
\begin{equation}
\label{eq_large_clusters}
\begin{split}
\sum_{C\in P(\WW_{\D\setminus \bar\L}(\o_{\partial_+\D})) \colon x\in \overline C}\hspace{-1cm}|\Phi_{\D\setminus \bar\L, \o_{\partial_+\D}}^p(C)|p^{-|\overline C|/(4d)}
&\leq \sum_{{\substack{C\in P(\WW_{\D\setminus \bar\L}(\o_{\partial_+\D}))  \colon x\in  \overline C}}} |\phi(C)|
\prod_{i = 1}^{n}p^{|W_i|/(2d)} \prod_{i = 1}^{n}p^{-|W_i|/(4d)} \\
&= \sum_{{\substack{C\in P(\WW_{\D\setminus \bar\L}(\o_{\partial_+\D}))  \colon x\in  \overline C}}} |\phi(C)|
\prod_{i = 1}^{n}p^{\frac{|W_i|}{2d}(1- \frac{1}{2})} \\
&= \sum_{{\substack{C\in P(\WW_{\D\setminus \bar\L}(\o_{\partial_+\D}))  \colon x\in  \overline C}}} |\phi(C)|
\prod_{i = 1}^{n}\sqrt{p}^{\frac{|W_i|}{2d}} \\
&=\sum_{C\in P(\WW_{\D\setminus \bar\L}(\o_{\partial_+\D}))  \colon x\in  \overline C} | \tilde \Phi_{\D\setminus \bar\L, \o_{\partial_+\D}}^{\sqrt{p}} (C)| \leq 1, 
\end{split}
\end{equation}
where $\tilde \Phi_{\D\setminus \bar\L, \o_{\partial_+\D}}^{\sqrt{p}} (C)$ denotes the potential of a cluster with weights given by
$$
\tilde z^{\sqrt{p}}_{W_k}(\o_{\partial_+\D}) = \sqrt{p}^{|W_k|/(2d)}.
$$
The first inequality in~\eqref{eq_large_clusters} is due to $|\overline C| \leq \sum_{i=1}^n|W_i|$. Moreover, assuming that $\sqrt{p}\le q_1$, the convergence criterion ~\cite[Theorem 5.4.]{friedli_velenik_2017} is satisfied for these adjusted weights (compare the argument in~\eqref{eq_s}). As a consequence, by~\cite[Statement (5.29)]{friedli_velenik_2017}, the contribution of all clusters containing a fixed site can be bounded by one, which gives the last inequality in~\eqref{eq_large_clusters}.
This gives that for any $R > 0$ and $p\le q^2_1$,
$$
\sum_{\substack{C\in P(\WW_{\D\setminus \bar\L}(\o_{\partial_+\D}))  \colon \\ x\in \overline C,\, | \overline C| \geq R}}|\Phi_{\D\setminus \bar\L, \o_{\partial_+\D}}^p(C)|p^{-R/(4d)}
\leq \sum_{\substack{C\in P(\WW_{\D\setminus \bar\L}(\o_{\partial_+\D}))  \colon \\ x\in \overline C,\, | \overline C| \geq R}}|\Phi_{\D\setminus \bar\L, \o_{\partial_+\D}}^p(C)|p^{-|\overline C|/(4d)}
\leq 1
$$
or equivalently
\begin{equation}
\label{eq_Bound_clusters}
\sum_{C\in P(\WW_{\D\setminus\bar\L}(\o_{\partial_+\D}))  \colon x\in \overline C,\, | \overline C| \geq R}|\Phi_{\D\setminus \bar\L, \o_{\partial_+\D}}^p(C)|
\leq p^{R/(4d)}.
\end{equation}
In words, we have achieved exponential suppression of the cluster potentials for large clusters that do not interact with the inner boundary.
\paragraph{Convergence of the fraction of partition functions:}
After having bounded the contribution of large clusters, we need to bound the non-suppressed polymers' weights in~\eqref{Eq_clus}. For this, note that for  $Q_k\in \QQ_{\D\setminus \bar\L}(\o_{\partial\L} \o_{\partial_+\D})$, we do not necessarily have $|z^p_{Q_k}(\o_{\partial\L} \o_{\partial_+\D})| \leq p^{|Q_k|/(2d)}$ like for the $W_k$, since the occupied sites leading to isolations may lie on $\partial_+\Lambda$, i.e., are part of the boundary condition. However, we have
\begin{equation}\label{eq_Q_k}
\begin{split}
|z^p_{Q_k}(\o_{\partial\L} \o_{\partial_+\D})|= \sum_{\o_{\D\setminus \bar\L}}\mu(\o_{\D\setminus \bar\L})\prod_{i\in Q_k}\one\{\o_{\D\setminus \bar\L}\in J_i(\o_{\partial\L} \o_{\partial_+\D})\}
\leq p^{|L_{Q_k} \cap \D\setminus \bar\L|} \leq p^{| L_{Q_k} |} p^{-|\partial_+\L \cap L_{Q_k} |},
\end{split}
\end{equation}
and since the $Q_k$ are pairwise disjoint and $\bigcup_{k=1}^n (\partial_+\L \cap L_{Q_k}) \subset  \partial_+\L$ and $|L_{Q_k}| \geq |Q_k|/(2d)$, we get
\begin{equation}\label{eq_Q_k}
\prod_{k=1}^n |z^p_{Q_k}(\o_{\partial\L} \o_{\partial_+\D})| \leq p^{\sum_{k=1}^n  |Q_k|/(2d)}p^{-|\partial_+\L|}.
\end{equation}
Next, let us denote by
$$
\C_{Q} := \{C \in P\big(\WW_{\D \sm \bar \L}(\o_{\partial_+ \D})\big) \colon C \cap \bigcup_{k=1}^n \overline Q_k \neq \emptyset\}
$$
the set of clusters interacting with the (non-suppressed) polymers $Q=\{Q_1,\dots,Q_n\}$. Then, 
$$
P\big(\WW_{\D\setminus\bar\L}(\o_{\partial_+\D})\big) \sm \C_{Q}  = P\big(\WW_{Q_1,\dots, Q_n, \D\setminus \bar\L}(\o_{\partial_+\D})\big),
$$
and we can now exploit cancellations by writing the fraction (\ref{Goal_1}) as
\begin{equation}\label{eq_poly_frac}
\begin{split}
&\tfrac{\sum_{\substack{Q_1, \dots, Q_n  \in \QQ_{\D\setminus \bar\L}(\o_{\partial\L}\o_{\partial_+\D})}}\prod_{k=1}^n z^p_{Q_k}(\o_{\partial\L} \o_{\partial_+\D})
\exp \big(\sum_{C \in P(\WW_{Q_1,\dots, Q_n, \D\setminus \bar\L}(\o_{\partial_+\D})) }\Phi^p(C) \big)}
{\sum_{\substack{Q_1, \dots, Q_n  \in \QQ_{\D\setminus \bar\L}(\tilde\o_{\partial\L}\o_{\partial_+\D})}}\prod_{k=1}^n z^p_{Q_k}(\tilde\o_{\partial\L} \o_{\partial_+\D})
\exp \left(\sum_{C \in P(\WW_{Q_1,\dots, Q_n, \D\setminus \bar\L}(\o_{\partial_+\D})) }\Phi^p(C) \right)} \\
&\qquad=\tfrac{\sum_{\substack{Q_1, \dots, Q_n  \in \QQ_{\D\setminus \bar\L}(\o_{\partial\L}\o_{\partial_+\D})}}\prod_{k=1}^n z^p_{Q_k}(\o_{\partial\L} \o_{\partial_+\D})
\exp \big(\sum_{C \in P(\WW_{\D\setminus\bar\L}(\o_{\partial_+\D})) }\Phi^p(C) - \sum_{C \in \C_{Q} }\Phi^p(C) \big)}
{\sum_{\substack{Q_1, \dots, Q_n  \in \QQ_{\D\setminus \bar\L}(\tilde\o_{\partial\L}\o_{\partial_+\D})}}\prod_{k=1}^n z^p_{Q_k}(\tilde\o_{\partial\L} \o_{\partial_+\D})
\exp \big(\sum_{C \in P(\WW_{\D\setminus\bar\L}(\o_{\partial_+\D})) }\Phi^p(C) - \sum_{C \in \C_{Q} }\Phi^p(C) \big)} \\
&\qquad=\tfrac{\sum_{\substack{Q_1, \dots, Q_n  \in \QQ_{\D\setminus \bar\L}(\o_{\partial\L}\o_{\partial_+\D})}}\prod_{k=1}^n z^p_{Q_k}(\o_{\partial\L} \o_{\partial_+\D})
\exp \big( - \sum_{C \in \C_{Q} }\Phi^p(C) \big)}
{\sum_{\substack{Q_1, \dots, Q_n  \in \QQ_{\D\setminus \bar\L}(\tilde\o_{\partial\L}\o_{\partial_+\D})}}\prod_{k=1}^n z^p_{Q_k}(\tilde\o_{\partial\L} \o_{\partial_+\D})
\exp \big( - \sum_{C \in \C_{Q} }\Phi^p(C) \big)}.
\end{split}
\end{equation}
For the sake of readability, we have omitted the indices of the cluster potentials. The next step is to verify convergence of the numerator and denominator respectively. For this, we can bound,
\begin{equation}\label{eq_C_1}
\sum_{C \in \C_{Q} } \lvert \Phi_{\D\setminus \bar\L, \o_{\partial_+\D}}^p(C) \rvert
\leq \bigl \lvert \bigcup_{k=1}^n \overline Q_k \bigr \rvert
\leq ((2d)^2 + 1)\sum_{k=1}^n |Q_k| =: ((2d)^2 + 1)m,
\end{equation}
where the absolute contribution of all clusters containing a fixed site can be bounded by one (cf.~the argument in~\eqref{eq_large_clusters}), and hence, the absolute contribution of clusters in $\C_{Q}$ can be bounded by $|\bigcup_{k=1}^n \overline Q_k|$, which is the first inequality in~\eqref{eq_C_1}.
Plugging in~\eqref{eq_Q_k} and~\eqref{eq_C_1}, we can bound the numerator of the last line in~\eqref{eq_poly_frac} by
\begin{equation}
\label{eq_Bound_numerator}
\begin{split}
p^{-| \partial_+\L|}\sum_{m = 0}^{\infty}\sum_{\substack{Q_1,\dots, Q_n \in \QQ_{\D\setminus \bar\L}(\tilde\o_{\partial\L}\o_{\partial_+\D}) \colon |\bigcup Q_k| = m}}p^{m/(2d)}
\e^{m((2d)^2 + 1)}.
\end{split}
\end{equation}
Given $m \in \N_{0}$, we need to count the number of sets of polymers $\{Q_1,\dots, Q_n\}$ with $\sum_k |Q_k| = m$ and such that for each $k$, $\overline Q_k \cap \partial_+ \L \neq \emptyset$. There are 
$$
{m + (|\partial_+\L|-1) \choose m} \leq \left( m + (|\partial_+\L| - 1) \right)^{(|\partial_+\L| - 1)}
$$
ways to distribute polymers of different sizes to disjoint starting nodes on $\partial_+\L$ such that the total size is given by $m$. For each such distribution with $n_1,\dots, n_{|\partial_+\L|} \in \N_{0}$, $\sum_{l=1}^{|\partial_+\L|}n_l = m$, there are at most $(9^d)^{m}$ different polymers (compare the argument following~\eqref{eq_s}). Therefore, we can bound~\eqref{eq_Bound_numerator} from above by
\begin{equation}
\label{eq_uniform_bound}
\begin{split}
p^{-|\partial_+\L|}\sum_{m = 0}^{\infty}\big( m + (|\partial_+\L| - 1) \big)^{(|\partial_+\L| - 1)} \left(p^{1/(2d)}
\e^{((2d)^2 + 1)} 9^d \right)^m,
\end{split}
\end{equation}
which is finite for sufficiently small $p\le q'_1$ independently of $\D$ and $S$.

\paragraph{Independence of boundary condition $\o_{\partial_+\D}$ as $\D\uparrow \Z^d$:}
We now analyze the fraction in (\ref{eq_poly_frac}) with respect to dependence on the outer boundary $\o_{\partial_+\D}$. 
For this, we further distinguish,
$$
\C_{Q, \cap \partial_+ \D} := \{C \in \C_{Q} \colon C \cap \partial_+ \D \neq \emptyset\},
$$
the subset of $\C_{Q}$ of polymers that also intersect the outer boundary and accordingly let
$$
\C_{Q, \not \cap \partial_+ \D} := \{C \in \C_{Q} \colon C \cap \partial_+ \D = \emptyset\}
$$ 
denote the clusters that do not reach the outer boundary. With this notation, we have
\begin{equation}\label{eq_cluster_decomp}
\C_{Q, \cap \partial_+ \D} \; \dot \cup \; \C_{Q, \not \cap \partial_+ \D} = \C_{Q},
\end{equation}
and, by the cluster decomposition~\eqref{eq_cluster_decomp}, the numerator in the last line of~\eqref{eq_poly_frac} takes the form

\begin{equation}\label{eq_clus_frac}
\begin{split}
\sum_{\substack{Q_1, \dots, Q_n \\ \text{pw.~comp.}}}\prod_{k=1}^n z^p_{Q_k}(\o_{\partial\L} \o_{\partial_+\D})
\exp \Big(- \sum_{\substack{C \in  \C_{Q, \cap \partial_+ \D} }}\Phi_{\D\setminus\bar\L,\o_{\partial_+\D}}^p(C)\Big)
\exp \Big(- \sum_{\substack{C \in \C_{Q, \not \cap \partial_+ \D}}}\Phi_{\D\setminus\bar\L}^p(C)\Big).
\end{split}
\end{equation}
Now, since clusters $C \in  \C_{Q, \cap \partial_+ \D} $ must suffice 
\begin{equation*}\label{eq_dist_bound}
\begin{split}
|\overline C| \geq \dist\Big(\bigcup_k \overline Q_k,\partial_+\D\Big)/4, 
\end{split}
\end{equation*}
the argument of the first exponential in~\eqref{eq_clus_frac} can be bounded from above by
\begin{equation}
\label{eq_bound_large_clusters}
\begin{split}
\sum_{\substack{C \in  \C_{Q, \cap \partial_+ \D}}} |\Phi_{\D\setminus\bar\L,\o_{\partial_+\D}}^p(C)| &\leq p^{\dist(\bigcup_k \bar Q_k, \partial_+\D)/(16d)},
\end{split}
\end{equation}
where we also used~\eqref{eq_Bound_clusters}.
This in particular implies that the first exponential term in~\eqref{eq_clus_frac} converges to one as $\D\uparrow\Z^d$. Moreover, using~\eqref{eq_Q_k}, also the contribution of $\prod_{k=1}^n z^p_{Q_k}(\o_{\partial\L} \o_{\partial_+\D})$ tends to zero, whenever it depends on $\o_{\partial_+\D}$. Hence, the function 
\begin{equation*}
\begin{split}
H_{\D,\o_{\partial_+\D}}(Q_1, \dots, Q_n)&:=\one\{Q_1, \dots, Q_n\subset\D\setminus\bar\L\}\prod_{k=1}^n z^p_{Q_k}(\o_{\partial\L} \o_{\partial_+\D})\\
&\qquad\times\exp \Big(- \sum_{\substack{C \in  \C_{Q, \cap \partial_+ \D} }}\Phi_{\D\setminus\bar\L,\o_{\partial_+\D}}^p(C)\Big)
\exp \Big(- \sum_{\substack{C \in \C_{Q, \not \cap \partial_+ \D}}}\Phi_{\D\setminus\bar\L}^p(C)\Big),
\end{split}
\end{equation*}
defined for any pairwise compatible set $\{Q_1, \dots, Q_n\}$ in $S$,
converges, as $\D\uparrow\Z^d$, to a function 
\begin{equation*}
\begin{split}
H(Q_1, \dots, Q_n):=\prod_{k=1}^n z^p_{Q_k}(\o_{\partial\L})\exp \Big(- \sum_{\substack{C \in \C_{Q}}}\Phi_{\bar\L^{\rm c}}^p(C)\Big),
\end{split}
\end{equation*}
independent of $\o_{\partial_+\D}$. But, since the bounds derived in~\eqref{eq_uniform_bound} are uniform in $\D$, we can employ the dominated-convergence theorem to conclude that 
\begin{equation*}
\begin{split}
\sum_{\substack{Q_1, \dots, Q_n \\ \text{pw.~comp.}}}H_{\D,\o_{\partial_+\D}}(Q_1, \dots, Q_n)\to\sum_{\substack{Q_1, \dots, Q_n \\ \text{pw.~comp.}}}H(Q_1, \dots, Q_n)\qquad\text{as }\D\uparrow\Z^d.
\end{split}
\end{equation*}
The same arguments also hold for the denominator in the last line of~\eqref{eq_poly_frac}, so we have finally arrived at our goal, the existence of the limit~\eqref{Goal_1} independent of the outer boundary condition.

\subsubsection{Existence of second-layer conditional expectations}
Finally, we wish to show that the limit as $\D\uparrow\Z^d$ of $\g^{S}_\D(F[\o'_\L]|\o_{\partial_+\D})$ exists and is independent of the boundary condition $\o_{\partial_+\D}$. 
For this, first note that, since the function $F[\o'_{\L}](\o)=F[\o'_{\L}](\o_{\partial\L})$ is local,
\begin{equation*}
\begin{split}
\g^S_\D(F[\o'_{\L}]| \o_{\partial_+\D})
= \sum_{\o_{\partial\L}} F[\o'_{\L}](\o_{\partial\L}) \g^S_\D(\o_{\partial\L}| \o_{\partial_+\D}),
\end{split}
\end{equation*}
and it suffices to consider $\g^S_\D(\o_{\partial\L}| \o_{\partial_+\D})$. 
But then, by the definition, we have
\begin{equation*}
\begin{split}
\g^S_\D(\o_{\partial \L} | \o_{\partial_+\D})
= \frac{Z_{\D\setminus \bar\L}(\o_{\partial\L} \o_{\partial_+\D})
}
{ \sum_{\ot_{\partial\L}} \mu(\ot_{\partial\L})Z_{\D\setminus \bar\L}(\ot_{\partial\L} \o_{\partial_+\D})
},
\end{split}
\end{equation*}
which does not depend on $\o_{\partial_+\D}$ as $\D \uparrow \Z^d$, by the previous step for $p\le p_1:=q^2_1\wedge q_1'$. This yields the result.

\bigskip
In particular, by the above, for all sufficiently small $p$, all $\L \Subset \zd$ and $\o'$, 
$$\g'(\o'_\L|\o'_{\L^c})=\lim_{\D\uparrow\Z^d}\g^{(\bar\o')^c}_\D(F[\o'_\L]|\o_{\D^c})$$ 
exists independently of $\o\in T^{-1}(\o')$. 

\subsubsection{Quasilocality of the specification}
For this, note that 
\begin{equation*}
\begin{split}
\sup_{\o',\, \eta' \colon \o'_{\D}  = \, \eta'_{ \D}}&|\g'(\o'_\L | \o'_{\L^c}) - \g'(\o'_\L | \eta'_{\L^c})|
\leq 2\sup_{\o',\, \o\in T^{-1}(\o')}|\g'(\o'_\L | \o'_{\L^c}) - \g'_{\o,\D}(\o'_\L | \o'_{\D \setminus \L})|\\
&+\sup_{\substack{\o',\,  \eta' \in \O' \colon \o'_{\D}  = \, \eta'_{\D}\\ \o\in T^{-1}(\o'),\, \eta\in T^{-1}(\eta') }}|\g^{(\bar\o')^{\rm c}}_\D(F[\o'_\L] | \o_{\partial_+\D}) - \g^{(\bar\eta')^{\rm c}}_{\D}(F[\o'_\L]  | \eta_{\partial_+\D})|,
\end{split}
\end{equation*}
where the first term on the right-hand side tends to zero as $\D$ tends to $\Z^d$ for sufficiently small $p$, by the cluster expansion arguments as presented above, which also gives the uniformity associated to $\o'$. For the second term on the right-hand side, 
it suffices to consider 
\begin{equation*}
\begin{split} 
|\g^{(\bar\o')^{\rm c}}_\D(\o_{\partial\L} | \o_{\partial_+\D}) - \g^{(\bar\eta')^{\rm c}}_{\D}(\o_{\partial\L} | \eta_{\partial_+\D})|,
\end{split}
\end{equation*}
which also becomes uniformly small in the boundary condition and the second-layer configurations, as $\D\uparrow\Z^d$ using the cluster-expansion arguments. This finishes the proof of Proposition~\ref{proposition_low}.

\subsection{Proof of Propositions~\ref{proposition_high} and~\ref{Prop_DoBound}}
\label{Sec_Proofs_High}
The distribution of a single site $i \in \zd$ depends only on finitely many spins, namely on all sites in the $l_1$-ball of $i$ of radius $2$ which we call the {\em dependence set} of site $i$. This dependence set being finite, we are dealing with a Markov field and we may apply Dobrushin-uniqueness techniques. However, as we see in the following example, the first-layer constraint model as formulated above does not directly permit applications of Dobrushin uniqueness or disagreement-percolation arguments. Indeed, writing $\eta^{\text{zero}}$ for the configuration without occupied sites, we have 
$$\g_0(0|\eta^{\text{zero}}) = 1\qquad\text{ and }\qquad\g_0(1|\eta^{\text{zero}}) = 0,$$
where we put $\g_0 := \g^{\Z^d}_{\{0\}}$ for the single-site specification kernel of the first-layer constraint model. On the other hand, denoting by $\eta^k$ the configuration that is fully empty except for a single occupied neighbor at $e_k$, the unit vector in the direction $k\in \{1, \dots, d\}$, we have
$$\g_0(0|\eta^k) = 0\qquad\text{ and }\qquad\g_0(1|\eta^k) = 1.$$
In particular, 
\begin{equation*}\begin{split}
\rho = \sup_{i \in \zd} &\max_{\eta,\tilde \eta\in \O}\Vert\g_i(\cdot|\eta)
-\g_i(\cdot|\tilde\eta)\Vert_{\text{TV}}=\frac{1}{2}\max_{\eta,\tilde \eta\in \O} (|\g_0( 0|\eta)-\g_0(0|\tilde\eta)| + |\g_0( 1|\eta)-\g_0(1|\tilde\eta)|)
=1,
\end{split}
\end{equation*}
independently of $p$, and hence the Dobrushin-uniqueness criterion or disagreement-percolation bounds cannot be satisfied.

However, we can rewrite in terms of a modified model, where this problem does not occur. Here, the idea is to form 2-by-1 pairs of sites that we think of as horizontal {\em dominos}, whose states we encode in terms of {\em pair-spin variables} $\xi_i$ with possible values $00,01,10,11$. We use the shorter single-digit notation for these pairs of symbols as $0,1,2,3$ in the sequel, simply reading them as two digits in a binary expansion. The new index set for the dominos is $\Z^{d-1} \times 2 \Z$, which is isomorphic to $\Z^d$. Let the axis alongside the dominos be denoted the {\em domino axis}. The first-layer constraint model is then equivalently described as a model on $\{0,1,2,3\}^{\Z^d}$ in terms of a translation-invariant single-site hardcore finite-range specification kernel $\varphi_i(\xi_i|\xi_{i^{\rm c}})$, where we suppress the dependence on the unfixed area.

\begin{figure}[t]
\centering
\includegraphics[width=.5\textwidth]{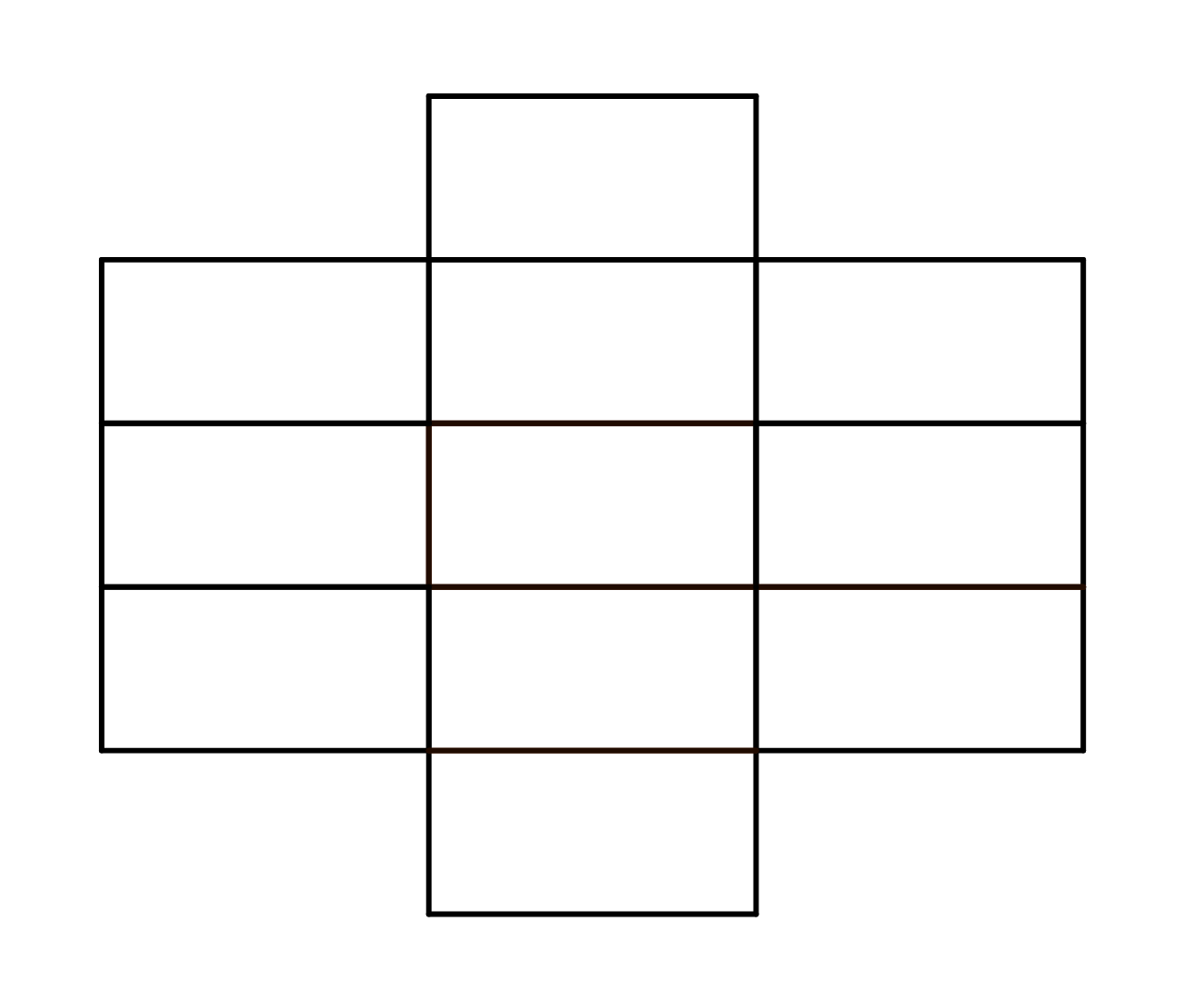}\hfill
\includegraphics[width=.5\textwidth]{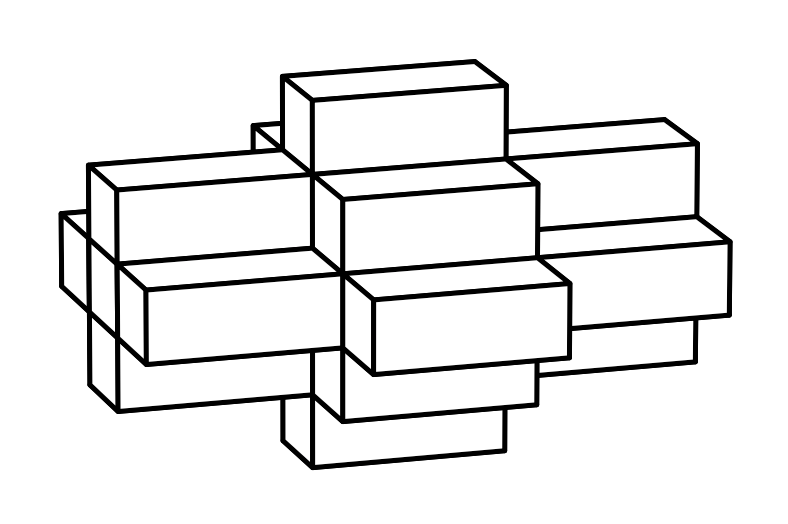}
\caption{The dependence sets $\protect V_0(2)$ (plus the center domino) and $\protect V_0(3)$ of the domino at the origin.}
\label{fig:figure_dependence2d}
\end{figure}

We first examine the dependence set $V_0(d)$ of $\varphi_0$, which now has a different form than in the introductory example. $V_0(d)$ contains the $2d$ adjacent dominos in its radius-$1$ $l_1$-boundary. Moreover, for each such domino, it contains the $2(d-1)$ adjacent dominos in all directions but the domino axis minus the center domino. Seeing the domino axis as the $x$-axis, we have a left and a right part which are symmetric and each contain $2(d-1) + 1$ dominos and a middle part containing $4(d-1) + 2(d-1)(d-2)$ dominos, see Figure~(\ref{fig:figure_dependence2d}). Altogether, we get
$$|V_0(d)| = 2(2(d-1) + 1) + 4(d-1) + 2(d-1)(d-2) = 2d^2 + 2d -2,$$ 
i.e., $V_0(2)$ contains $10$ sites of dominos and $V_0(3)$ contains $22$ sites.

An important observation is that the single-site kernels $\varphi_i(\cdot|\xi)$ are completely specified by asking which values from $\{0,1,2,3\}$ are assigned to non-zero probabilities. For the equivalence class of boundary conditions, for which the values $0$ and $3$ are allowed, while $1$ and $2$ are forbidden as they describe isolated particles, we write $+,-,-,+$. For example, the all-zero boundary condition $\xi^{\text{zero}}$ belongs to the class $+,-,-,+$, since
$$\varphi_0(1|\xi^{\text{zero}})=\varphi_0(2| \xi^{\text{zero}})=0$$
as isolated occupied sites are forbidden, while 
$$\varphi_0(0| \xi^{\text{zero}})=\frac{(1-p)^2}{p^2+(1-p)^2}\qquad\text{ and }\qquad\varphi_0(3| \xi^{\text{zero}})=\frac{p^2}{p^2+(1-p)^2}$$
are determined as the Bernoulli measure conditioned on the allowed values $0$ and $3$. 
The corresponding probability vector for the values 
$0,1,2,3$ takes the form 
\begin{equation*}\begin{split}
\varphi_0(\cdot|\xi^{\text{zero}})=\frac{1}{p^2+(1-p)^2}((1-p)^2,0,0,p^2).
\end{split}
\end{equation*}
To compare, for the fully occupied boundary $\xi^{\text{one}}$, all pairs inside are allowed, 
which is why it belongs to the class $+,+,+,+$, 
and we have 
\begin{equation*}\begin{split}
\varphi_0(\cdot |\xi^{\text{one}})=((1-p)^2,p(1-p),p(1-p),p^2).
\end{split}
\end{equation*}
We note that the value $3$ is allowed for all possible boundary conditions, as it prevents isolations both in the center and on the boundary. This implies that for all possible boundary conditions $\xi$, we have that
\begin{equation*}\begin{split}
\varphi_0(\cdot|\xi) \rightarrow (0,0,0,1),
\end{split}
\end{equation*}
as $p \uparrow 1$. This means that single-domino conditional measures become concentrated and the system enters a strong-field regime. Indeed, this feature also makes the essential quantities~\eqref{def_DM} and \eqref{def_rho} in the Dobrushin approach and the percolation framework decrease for large values of $p$. Note further that the last observation means that all probability vectors must be in classes of the form $\cdot, \cdot, \cdot, +$. However, there is no boundary condition belonging to the string $-,+,+,+$, independent of the dimension $d$. Indeed, if domino value $0$ leads to an isolation on the boundary, this isolation will have to occur either for domino values $1$ or $2$ contradicting the fact that both of them yield a non-zero probability for the string $-,+,+,+$. The remaining $7$ strings can all occur.

Thinking of the lattice as a dependence graph $G = (\zd, E_{\text{dep}})$, where each $i$ is connected to all $j \in V_i(d)$, we are dealing with a Markov field and would like to apply uniqueness criteria.

\subsubsection{Uniqueness via the Dobrushin-uniqueness criterion}\label{Sec_Dob1}
\begin{figure}[t]
\centering
\includegraphics[width=1\textwidth]{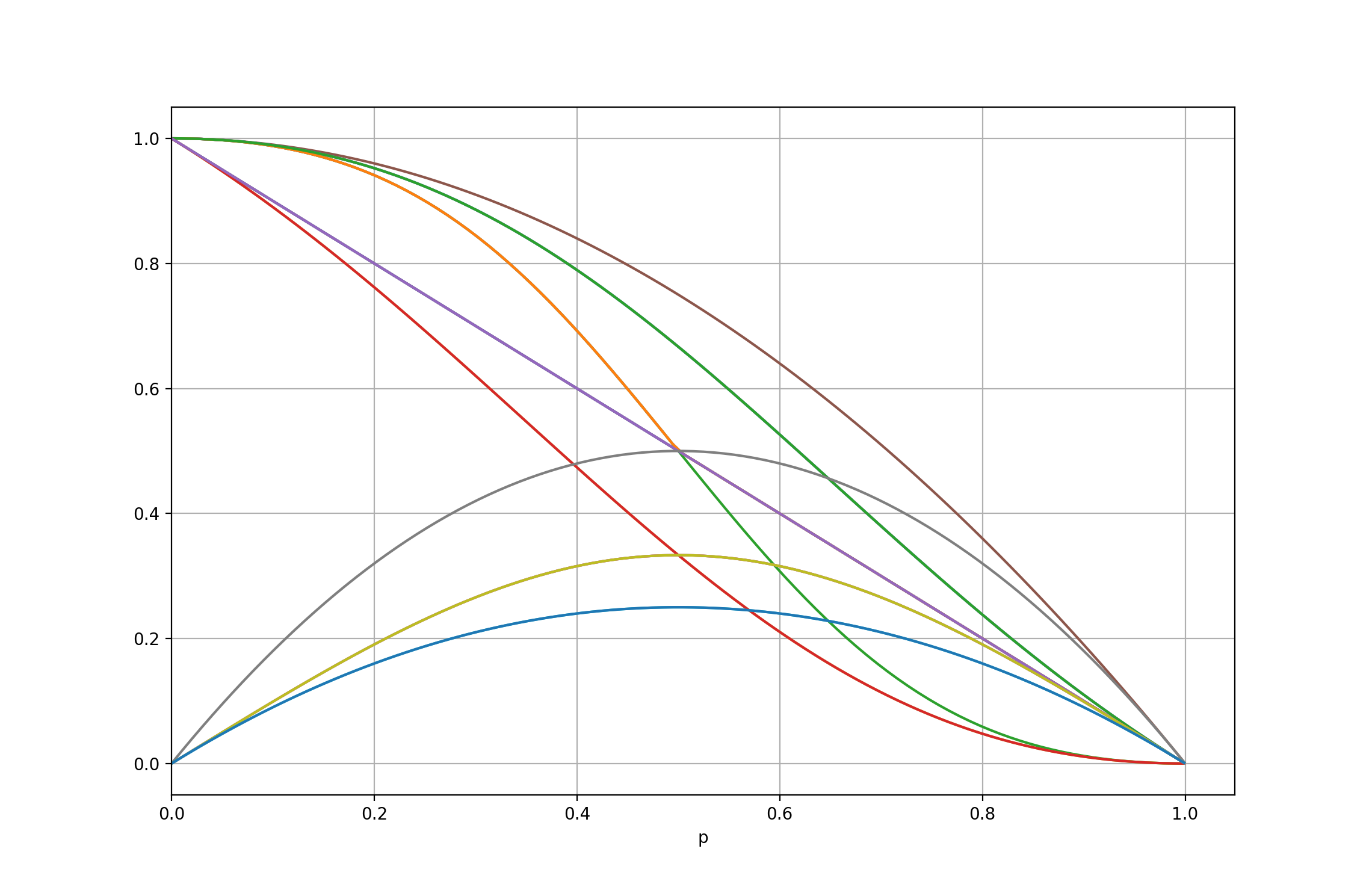}\hfill
\caption{All possible total variational distances. The $\protect {7 \choose 2} = 21$ possible pairs of strings yield $\protect 8$ distinct curves. The $\protect 4$ top ones in the high-density regime (i.e., for $\protect  p > 0.6479$) are given by rational functions $\protect \rho \geq q  \geq u \geq v$ (red, grey, green, purple) in (\ref{eq_polynomials}). The code for this computation can be found in~\cite{dominos}.}
\label{fig:figure2}
\end{figure}

Let us reintroduce the unfixed area $S\subset\Z^d$. We determine the Dobrushin matrix
\begin{equation}\label{def_DM}
C^S_{i,j}(p, d) := \max_{\xi,\tilde\xi\in\{0,1\}^S\colon \xi_{j^c} = \tilde\xi_{j^c}}\|\varphi^S_i(\cdot|\xi)
-\varphi^S_i(\cdot|\tilde\xi)\|_{\text{TV}}
\end{equation}
for the domino specification in order to compute the Dobrushin constant
\begin{equation}\label{def_DC}
c^S(p, d):=\sup_{i\in S}\sum_{j \in V_i(d)}C^S_{i,j}(p, d).
\end{equation}
Then, we are interested in the lowest threshold $p^{\rm d}_{\rm c}(d)$, such that $c(p, d):=\sup_{S\subset\Z^d} c^S(p, d)< 1$ for $p \geq p^{\rm d}_{\rm c}(d)$. Note that this then corresponds to a uniform bound in $S$. A simple but non-optimal bound is given by $c(p, d)\leq |V_0(d)| \rho(p)$, where 
\begin{equation}\label{def_rho}
\begin{split}
\rho(p):=\max_{\xi,\tilde \xi\in\O } \|\varphi^{\Z^d}_i(\cdot|\xi)
-\varphi^{\Z^d}_i(\cdot|\tilde \xi)\|_{\text{TV}}
=\frac{1}{2}\max_{\xi,\tilde \xi\in \O}\sum_{a=0,1,2,3}|\varphi^{\Z^d}_i(a|\xi),
-\varphi^{\Z^d}_i(a|\tilde \xi)|,
\end{split}
\end{equation}
which is independent of $i \in \zd$ due to translation invariance of the kernels. By the above, we may determine $\rho(p)$ by computing the total variational distances of at most ${7 \choose 2} = 21$ pairs of probability vectors. A straightforward computation shows that for all $p \in [0,1]$, the maximum is attained by the probability vectors corresponding to the strings $-,-,-,+$ and $+,+,+,+$, see Figure~\ref{fig:figure2}, 
and thus
\begin{equation}
\begin{split}
\rho(p) 
= & \, \|((1-p)^2,p(1-p), p(1-p), p^2) - (0,0,0,1)\|_{\text{TV}}=1 - p^2.
\end{split}
\end{equation}
In particular, for $p\ge p_{\rm c}(d)$, with 
$$
p_{\rm c}(d):=\sqrt{\frac{2d^2 + 2d - 2 }{2d^2 + 2d - 3}},
$$ 
we are in the Dobrushin-uniqueness regime, see~\cite[Theorem 8.7 and Equation 8.25]{Ge11}, and thus, the unique existence of the
$\lim_{\D\uparrow\Z^d}\varphi^{(\bar\o')^c}_\D(F[\o'_\L]|\o_{\D^c})$ is guaranteed independently of $\o'$ and $\o_{\D^c}$.  

\subsubsection{Quasilocality of the specification}
What remains to be done in order to finish the proof of Proposition~\ref{proposition_high} is to establish quasilocality for the specification $\g'$. For this, we let $s$ denote the $\ell_\infty$ metric on $\Z^d$ and define $s(\L,\D)=\inf\{s(i,j)\colon i\in \L, j\in \D\}$. Then we have the following result, which is equivalent to the corresponding result in the companion paper~\cite[Lemma 3.5]{jahnel2021gibbsianness}.
\begin{lem}\label{lem_Quasi}
For $p> p_{\rm c}(d)$ there exist constants $C,c>0$ such that for all $\L\subset\D\Subset\Z^d$ and all configurations $\o'$ and $\eta'$ with $\o'_\D=\eta'_\D$ we have that 
\begin{equation*}
|\g'_\L(\o'_\L|\o'_{\L^c})-\g'_\L(\o'_\L|\eta'_{\L^c})|\leq  C |\Lambda| e^{-c s(\Lambda, \D^c)}.
\end{equation*}
In particular, the specification $\g'$ is quasilocal. 
\end{lem}
We briefly sketch the proof here for completeness. 
\begin{proof}
Note that in the regime $p> p_{\rm c}(d)$, we can represent $\g'_\L(\o'_\L|\o'_{\L^c})$ via the unique infinite-volume Gibbs measure $\mu^{(\bar\o')^{\rm c}}(F[\o'_\L])$. Now, using the criterion~\cite[Remark 8.26]{Ge11} applied to~\cite[Theorem 8.20]{Ge11}, we have that 
\begin{equation*}
|\mu^{(\bar\o')^{\rm c}}(F[\o'_\L])-\mu^{(\bar\eta')^{\rm c}}(F[\o'_\L])|\le D(\L,\D), 
\end{equation*}
where $D(\L,\D)=\sum_{i\in \L, j\in \D^c}\big(\sum_{n\ge 0}C^n\big)_{i,j}$ with $C^n$ the $n$-th power of the Dobrushin matrix $C = (C_{i,j}(p, d))_{i,j\in \Z^d}$ as defined in~\eqref{def_DM}.  Now choose $c>0$ sufficiently small such that $p>e^cp_{\rm c}(d)$, then, by~\cite[Remark 8.26]{Ge11}, 
\begin{equation*}
D(\L,\D)\le C|\L|\e^{-c d(\L,\D^c)}
\end{equation*}
for some finite $C>0$ and the proof is finished. 
\end{proof}

\subsubsection{Proof of Proposition~\ref{Prop_DoBound}}
Recall the definition of the Dobrushin constant $c(p,d)$ from~\eqref{def_DC}. Then, the statement of Proposition~\ref{Prop_DoBound} follows directly from the statement of the following lemma. 
\begin{lem} We have that 
\begin{equation*}
\begin{split}
c(p,d) &= (1-p^2)(2(d-1)(d-2)) + 4(d-1)p(1-p) + 2\frac{1-p}{1 - p(1-p)}+ (1-p)(6(d-1)).
\end{split}
\end{equation*}
\end{lem}
\begin{proof}
Without loss of generality, let the domino axis point along the first unit vector $e_1$. Due to symmetries, for many dominos $j \in V_0(d)$, the contributions $C_{0,j}(p)$ are the same. Therefore, we divide $V_0(d)$ into the $5$ disjoint classes, see Figure~\ref{fig:figure_dependence2d}, 
\begin{enumerate}
\item[] $V_1(d) := \{\pm e_i \colon i = 2, \dots, d\}$, $|V_1(d)| = 2(d-1)$, the direct neighbors of the center in the non-domino directions,
\item[] $V_2(d) := \{\pm e_1\}$, $|V_2(d)| = 2$, the direct neighbors of the center in the domino directions, 
\item[] $V_3(d) := \{\pm e_1 \pm e_i \colon  i = 2, \dots, d\}$, $|V_3(d)| = 4(d-1)$, the direct neighbors in the non-domino directions of sites in $V_2(d)$, 
\item[] $V_4(d) := \{\pm 2 e_i \colon i = 2, \dots, d\}$, $|V_4(d)| = 2(d-1)$, the distance-$2$ sites of the center in the non-domino directions, and 
\item[] $V_5(d) := \{\pm e_i \pm e_j \colon i, j = 2, \dots, d\} \setminus V_4$, $|V_5(d)| = 2(d-1)(d-2)$, the direct neighbors in the non-domino directions of $V_1$ inside the $l^{\infty}$-ball of radius $2$.
\end{enumerate}
Next, we define the rational functions
\begin{equation}\label{eq_polynomials}
\begin{split}
\rho(p) &:= 1 - p^2 = \text{TV}((-,-,-,+), (+,+,+,+)), \\
q(p) &:= 2p(1-p)= \text{TV}((+,-,-,+), (+,+,+,+)), \\
u(p) &:= \frac{1-p}{1 - p(1-p)}=  \text{TV}((-,-,+,+), (+,+,-,+))= \text{TV}((-,+,-,+), (+,-,+,+))\\
&=\text{TV}((-,-,-,+), (+,-,+,+))=  \text{TV}((-,-,-,+), (+,+,-,+)), \\
v(p) &:= 1-p=  \text{TV}((-,-,-,+), (-,-,+,+)) =  \text{TV}((-,-,-,+), (-,+,-,+)) \\
&=  \text{TV}((-,-,+,+), (-,+,-,+)) =  \text{TV}((-,-,+,+), (+,-,-,+)) \\
&=  \text{TV}((-,-,+,+), (+,+,+,+)) =  \text{TV}((-,+,-,+), (+,-,-,+)) \\
&=  \text{TV}((-,+,-,+), (+,+,+,+)), 
\end{split}
\end{equation}
where $\text{TV}(a,b)=\sum_{i=1,\dots,4}|a_i-b_i|/2$ denotes the total-variational distance of the probability vectors $a,b$ belonging to the classes defined by the strings. For example, we have that 
\begin{equation*}
\begin{split}
\text{TV}((-,-,-,+), (+,+,+,+)) &= \text{TV}((0,0,0,1), ((1-p)^2,p(1-p),p(1-p),p^2)) \\ 
&=\frac{1}{2}\left((1-p)^2 + 2p(1-p) + 1-p^2 \right) = 1-p^2.
\end{split}
\end{equation*}
Note that, for $p > 0.6479$, we have $\rho(p) > q(p) > u(p) > v(p)$, see Figure~\ref{fig:figure2}. 
In the following, speaking of sites refers to indices of the original lattice, i.e., we say each domino has a left and a right site. Speaking of dominos refers to indices in the domino lattice. Imagine the hyperplane in the original lattice orthogonal to the domino axis and separating the center domino's two sites. The left {\em halfspace or side} denotes sites on the left side of this hyperplane (towards the negative domino axis) while the right halfspace (side) refers to sites on the right side of this hyperplane.

The strategy is to go through the polynomials (\ref{eq_polynomials}) in decreasing order and to try to construct boundary conditions such that, if subjected to a domino flip in some $V_i, i=1, \dots, 5$, the total-variational distance of that flip is given by the polynomial.

We start by $\rho$ and check if there is a boundary condition $\xi$ of class $-,-,-,+$, such that a single domino flip yields class $+,+,+,+$. Class $-,-,-,+$ implies that in each halfspace there is an isolation on the boundary, while the isolations of course occur on different dominos in $V_1(d) \cup V_2(d)$. For $+,+,+,+$ there cannot be isolations on the boundary. A domino whose flipping transforms $-,-,-,+$ into $+,+,+,+$ needs to be adjacent to both dominos creating the isolations. This excludes dominos from $V_1(d), V_2(d), V_3(d)$ and $V_4(d)$. However, by flipping the boundary $\xi$ consisting only of zeros, except for values $1$ and $2$ at dominos $e_2$ and $e_3$ at domino $e_2 + e_3 \in V_5(d)$, from $0$ to $3$ leads to a change from $-,-,-,+$ to $+,+,+,+$. Of course, by the symmetry of $V_5(d)$, one can construct analogous boundaries for flips at all other dominos of $V_5(d)$.

Checking for $q$, we find that a boundary from $+,-,-,+$ means that there are no isolations on the boundary, but center domino values $1$ and $2$ create isolations on the domino, meaning that there are only unoccupied sites around each of the center domino sites. Transformation to $+,+,+,+$ thus requires a domino from $V_1(d)$ to be flipped from $0$ to $3$ and flips on $V_2(d)$, $V_3(d)$, $V_4(d)$ and $V_5(d)$ cannot lead to this transformation. Note that for $d=2$, $V_5(d) = \emptyset$.

Checking for $u$ yields the most technical case. We start by flipping $-,-,+,+$ to $+,+,-,+$. $-,-,+,+$ means that there is at least one isolation on the left-side boundary and no isolation on the right side boundary. $+,+,-,+$ implies that there is no isolation on the boundary, at least one non-isolated site on the right boundary and only unoccupied sites around the left domino site. A change from one to the other requires the domino with the isolated site on the left boundary to be flipped. This excludes dominos from $V_3(d), V_4(d), V_5(d)$, but flips in $V_1(d)$ and $V_2(d)$ are possible. An  example boundary is given by $\xi$ in $-,-,+,+$ with value zero everywhere except for a value $1$ at position $-e_1$ and a value $3$ at position $e_1$, flipped at $-e_1 \in V_2(d)$ from value $1$ to value $0$. The transformation $-,+,-,+$ to $+,-,+,+$ is symmetric with respect to the hyperplane orthogonal to $e_1$ and thus only yields domino flips at $V_1(d)$ or $V_2(d)$, too.

We continue with class $-,-,-,+$ which demands that there is at least one isolation on the boundary on each side. For $+,-,+,+$, there cannot be isolation on the boundary, the right domino site is surrounded by unoccupied sites, while on the left side, there is at least one pair of neighboring occupied sites on the boundary. To transform, two things need to happen: First, the isolation on the left side needs to get another occupied neighbor, which demands a domino flip in $V_3(d) \cup V_4(d) \cup V_5(d)$. Second, the domino containing the isolation on the right side which is thus of value $1$ or $2$ and belongs to $V_1(d) \cup V_2(d)$ needs to be flipped to $0$. It is impossible for both to happen with a single domino flip. Due to symmetry, the same holds for the the transformation $-,-,-,+$ to $+,+,-,+$.

Finally, checking for $v$ yields that we do not need to go through all possible transformations. It is enough to show that there exist boundary conditions such that a domino flip in $V_3(d)$ and $V_4(d)$ respectively yields a variational distance of weight $v$. For a flip at $e_1 + e_i \in V_3(d)$, $i = 2, \dots, d$, consider the change $-,-,-,+$ to $-,-,+,+$ by the boundary $\o$ with zeros except for value $1$ at $-e_1$, value $2$ at $e_1$. Flipping at $e_1 + e_i$ from $0$ to $2$ yields the desired transformation. The same construction can be done on the left side applying the transformation $-,-,-,+$ to $-,+,-,+$ to reach of the remaining dominos in $V_3(d)$.

Consider again the change $-,-,-,+$ to $-,-,+,+$, this time for the boundary $\xi$ consisting of zeros except for value $1$ at $e_i$ and value $2$ at $-e_i$ for some $i \in 2, \dots, d$. Flipping the domino at $2e_i \in V_4(d)$ from value $0$ to value $1$ serves the purpose. Again, the left side dominos of $V_4(d)$ can be reached by applying the same construction in the symmetric case $-,-,-,+$ to $-,-,+,+$.

Altogether, we get 
\begin{equation*}
\begin{split}
c(p, d) &=\rho(p) |V_5(d)| + q(p) |V_1(d)| + u(p) |V_2(d)| + v(p) (|V_3(d)| + |V_4(d)|) \\
&=(1-p^2)(2(d-1)(d-2)) + 4(d-1)p(1-p) + 2\frac{1-p}{1 - p(1-p)}+ (1-p)(6(d-1))
\end{split}
\end{equation*}
and this finishes the proof. 
\end{proof}

\subsection{Proof of Proposition~\ref{Prop_DisBound}}\label{Sec_Prop_DisBound}
The proof is based on disagreement-percolation bounds for general graphs. 
\begin{proof}
By the main result of~\cite{BeMa94} and the site-percolation bound 
\begin{equation}\label{eq_percolation_bound}
p_c \geq \frac{1}{\sup_i |N_i| - 1},
\end{equation}
where $N_i := \partial_+\{i\}$ denotes the set of neighbors of $i$, applied to the locally finite dependence graph induced by the domino model, we have uniqueness once $\rho(p) < \frac{1}{|V_0(d)| - 1} \leq p_c$, where $\rho(p)$ is defined in~\eqref{def_rho}. Then, plugging in $\rho(p) = 1-p^2$ and $|V_0(d)| = 2d^2 + 2d - 2$ yields the unique existence of infinite-volume Gibbs measure for the domino model. By the equivalence between the domino model and the first-layer constraint model, 
absence of a phase transition also leads to convergence of
$\g_\D^S(F[\o'_{\L}]|\o_{\partial_+ \D})$ as $\D\uparrow\Z^d$, independent of the boundary condition $\o$ and the unfixed area $S$. 
\end{proof}

\section*{Acknowledgements}
This work was funded by the German Research Foundation under Germany's Excellence Strategy MATH+: The Berlin Mathematics Research Center, EXC-2046/1 project ID: 390685689 and the German Leibniz Association via the Leibniz Competition 2020.

\bibliography{Jahnel}
\bibliographystyle{abbrv}
\end{document}